\documentclass{agtart_a}
\pdfoutput=1

\usepackage{graphicx}

\title[Conjugacy of $2$--spherical subgroups of Coxeter groups and 
parallel walls]{Conjugacy of $2$--spherical subgroups of Coxeter groups\\ 
and parallel walls}

\author{Pierre-Emmanuel Caprace}
\givenname{Pierre-Emmanuel}
\surname{Caprace}
\address{D\'epartement de Math\'ematiques\\
Universit\'e Libre de Bruxelles, CP216\\\newline
Bd du Triomphe\\1050 Bruxelles\\Belgium}
\email{pcaprace@ulb.ac.be}
\urladdr{}

\volumenumber{6}
\issuenumber{}
\publicationyear{2006}
\papernumber{70}
\startpage{1987}
\endpage{2029}

\doi{}
\MR{}
\Zbl{}


\received{2 August 2005}
\revised{31 August 2006}
\accepted{4 October 2006}
\published{14 November 2006}
\publishedonline{14 November 2006}
\proposed{}
\seconded{}
\corresponding{}
\editor{CPR}
\version{}

\arxivreference{math.GR/0508057}

\let\xysavmatrix\xymatrix
\def\xymatrix{\disablesubscriptcorrection\xysavmatrix}
\AtBeginDocument{\let\bar\wbar\let\tilde\wtilde\let\hat\what}

\makeatletter
\def\dnewtheorem#1[#2]#3{\newtheorem{#1}{#3}
\expandafter\let\csname c@#1\endcsname\c@thmintro}
\def\cnewtheorem#1[#2]#3{\newtheorem{#1}{#3}
\expandafter\let\csname c@#1\endcsname\c@prop}
\makeatother

\theoremstyle{plain}
\newtheorem{prop}{Proposition}
\cnewtheorem{thm}[prop]{Theorem}
\cnewtheorem{cor}[prop]{Corollary}
\cnewtheorem{lem}[prop]{Lemma}
\newtheorem{thmintro}{Theorem}

\makeautorefname{thmintro}{Theorem}
\dnewtheorem{corintro}[thmintro]{Corollary}

\makeautorefname{corintro}{Corollary}

\theoremstyle{definition}

\theoremstyle{remark}

\newtheorem*{remark}{Remark}
\newtheorem{remarknb}{Remark}

\newcommand{\Hyp}{\mathbb{H}}

\newcommand{\N}{\mathbb{N}}
\newcommand{\GL}{\mathrm{GL}}

\newcommand{\la}{\langle}
\newcommand{\ra}{\rangle}
\newcommand{\inv}{^{-1}}

\newcommand{\CAT}[1]{\mathsf{CAT}(#1)}%
\def\bs#1.{
              \def\temp{#1}
              \ifx\temp\empty
                   \mathcal{B}
              \else
                   \mathcal{B}(#1)
              \fi
}
\DeclareMathOperator{\pr}{pr} 
\DeclareMathOperator{\Stab}{Stab}

\DeclareMathOperator{\Pc}{Pc}

\DeclareMathOperator{\Zrk}{\mathbb{Z}-rk}
\DeclareMathOperator{\Isom}{Isom}

\def\min{\mathop{\mathrm{min}}\nolimits}

\def\max{\mathop{\mathrm{max}}\nolimits}

\begin{document}

\begin{asciiabstract}
Let (W,S) be a Coxeter system of finite rank (ie |S| is finite) and
let A be the associated Coxeter (or Davis) complex. We study chains of
pairwise parallel walls in A using Tits' bilinear form associated to
the standard root system of (W,S). As an application, we prove the
strong parallel wall conjecture of G Niblo and L Reeves [J Group
Theory 6 (2003) 399--413]. This allows to prove finiteness of the
number of conjugacy classes of certain one-ended subgroups of W,
which yields in turn the determination of all co-Hopfian Coxeter
groups of 2--spherical type.
\end{asciiabstract}

\begin{htmlabstract}
Let (W,S) be a Coxeter system of finite rank (ie &#x2223;S&#x2223; is finite)
and let A be the associated Coxeter (or Davis) complex. We
study chains of pairwise parallel walls in A using Tits'
bilinear form associated to the standard root system of (W,S). As an
application, we prove the strong parallel wall conjecture of G&nbsp;Niblo
and L&nbsp;Reeves [J Group Theory 6 (2003) 399&ndash;413]. This allows to prove
finiteness of the number of conjugacy classes of certain one-ended
subgroups of W, which yields in turn the determination of all
co-Hopfian Coxeter groups of 2&ndash;spherical type.
\end{htmlabstract}

\begin{webabstract}
Let $(W,S)$ be a Coxeter system of finite rank (ie $|S|$ is finite)
and let $\mathcal{A}$ be the associated Coxeter (or Davis) complex. We
study chains of pairwise parallel walls in $\mathcal{A}$ using Tits'
bilinear form associated to the standard root system of $(W,S)$. As an
application, we prove the strong parallel wall conjecture of G~Niblo
and L~Reeves [J Group Theory 6 (2003) 399--413]. This allows to prove
finiteness of the number of conjugacy classes of certain one-ended
subgroups of $W$, which yields in turn the determination of all
co-Hopfian Coxeter groups of $2$--spherical type.
\end{webabstract}

\begin{abstract}
Let $(W,S)$ be a Coxeter system of finite rank (ie $|S|$ is finite)
and let $\mathscr{A}$ be the associated Coxeter (or Davis) complex. We
study chains of pairwise parallel walls in $\mathscr{A}$ using Tits'
bilinear form associated to the standard root system of $(W,S)$. As an
application, we prove the strong parallel wall conjecture of G~Niblo
and L~Reeves \cite{NR03}. This allows to prove finiteness of the
number of conjugacy classes of certain one-ended subgroups of $W$,
which yields in turn the determination of all co-Hopfian Coxeter
groups of $2$--spherical type.
\end{abstract}

\maketitle

\section{Introduction}

\subsection{Conjugacy of $2$--spherical subgroups}\label{sect:2sph}

A group $\Gamma$ is called \textit{$2$--spherical} if it possesses a
finite generating set $T$ such that any pair of elements of $T$
generates a finite subgroup. By Serre \cite[Section 6.5,
Corollaire~2]{SerreArbres}, a $2$--spherical group enjoys property
(FA); in particular, it follows from Stalling's theorem that it is
one-ended.

In the literature, a Coxeter group $W$ is called $2$--spherical if it
has a Coxeter generating set $S$ with the property that any pair of
elements of $S$ generates a finite subgroup. If $W$ has a Coxeter
generating set $S$ such that some pair of elements of $S$ generates an
infinite subgroup, then it is easy to see that $W$ splits
non-trivially as an amalgamated product of standard parabolic
subgroups, and hence $W$ does not have Serre's property (FA). This
shows that for Coxeter groups, the usual notion of $2$--sphericity
coincides with the notion introduced above.

A theorem stated by M~Gromov \cite[5.3, C$'$]{Gromov} and proved by
E~Rips and Z~Sela \cite[Theorem~7.1]{RS94} asserts that a hyperbolic
group contains only finitely many conjugacy classes of subgroups
isomorphic to a given finitely generated torsion-free one-ended group
$\Lambda$. Another related result, stated in \cite{Gromov} and proved
by T~Delzant \cite{Delzant91}, is that a torsion-free hyperbolic group
$\Gamma$ has finitely many conjugacy classes of one-ended
two-generated subgroups. Furthermore, it was shown by I~Kapovich and
R~Weidmann \cite[Corollary~1.5]{KapovichWeidmann04} that if $\Gamma$
is locally quasiconvex, then $\Gamma$ has finitely many conjugacy
classes of one-ended $l$--generated subgroups for each $l \geq
1$. Results of this kind corroborate somehow the abundance of free
subgroups in hyperbolic groups.  Comparing Coxeter groups, the
abundance of free subgroups is also an established fact (more
precisely: Coxeter groups are either virtually abelian or virtually
map onto non-abelian free groups by Margulis and Vinberg
\cite[Corollary~2]{MarVin00}) and it is desirable to complement this
fact with a precise finiteness property. In this direction, we obtain
the following:

\begin{thmintro}\label{thm:FiniteConjug:2sph}$\phantom{9}$
\begin{itemize}
\item[(i)] $W$ contains finitely many conjugacy classes of $2$--spherical subgroups if and only if $W$ has no
parabolic subgroup of irreducible affine type and rank~$\geq 3$.

\item[(ii)] $W$ contains finitely many conjugacy classes of $2$--spherical reflection subgroups $\Gamma$ such
that $\Gamma$ has no nontrivial free abelian normal subgroup (equivalently: $\Gamma$ has no direct component of
affine type).

\item[(iii)]  $W$ contains finitely many conjugacy classes of $2$--spherical subgroups $\Gamma$ such that
$\Gamma$ has no nontrivial free abelian normal subgroup and $\Zrk(\Gamma) \geq \Zrk(W) -1$.

\item[(iv)] Assume that the centralizer of every irreducible affine parabolic subgroup of rank~$\geq 3$ is
finite. Then $W$ contains finitely many conjugacy classes of $2$--spherical subgroups $\Gamma$ such that $\Gamma$
is not infinite virtually abelian.
\end{itemize}
\end{thmintro}

By a reflection subgroup, we mean a subgroup generated by reflections; see \fullref{sect:ReflSubgps} below for
more information on reflection subgroups. Here $\Zrk(\Gamma)$ denotes the maximal rank of a free abelian
subgroup of $\Gamma$. Note that the condition on the centralizer of an affine parabolic subgroup in (iv) is
immediate to check in the Coxeter diagram of $(W,S)$ (see \fullref{lem:Normalizer} below): the condition holds
whenever for each irreducible affine subdiagram $D$ of rank~$\geq 3$, the subdiagram induced on the set
$D^\perp$ consisting of those vertices which are not connected to $D$ is of spherical type.

The existence of homotheties in Euclidean space implies that any Coxeter group of affine type is isomorphic to a
proper reflection subgroup. From this fact, it is easy to deduce that if $W$ has a reflection subgroup $W_0$ of
affine type, then $W$ has infinitely many conjugacy classes of reflection subgroups of the same type as $W_0$.
This explains why the number of conjugacy classes of $2$--spherical subgroups of $W$ necessarily depends on
affine parabolic subgroups of $W$, as it appears in \fullref{thm:FiniteConjug:2sph}. However, affine
parabolic subgroups are not the only source of free abelian subgroups in Coxeter groups; in particular, many
non-hyperbolic Coxeter groups possess no parabolic subgroup of irreducible affine type and rank~$\geq 3$.

A group which is not isomorphic to any of its proper subgroups
(resp.\ quotients) is called \textit{co-Hopfian}
(resp.\ \textit{Hopfian}). By a theorem of Mal'cev, any finitely
generated residually finite group is Hopfian; hence Coxeter groups are
Hopfian. By Rips and Zela \cite[Theorem~3.1]{RS94}, rigid hyperbolic
groups are co-Hopfian. We have just seen that affine Coxeter groups
are not co-Hopfian. As a consequence of
\fullref{thm:FiniteConjug:2sph}, one obtains the following:

\begin{corintro}\label{cor:Hopfian}
Suppose that $W$ is $2$--spherical. Then $W$ is co-Hopfian if and only if $W$ has no nontrivial free
abelian normal subgroup or, equivalently, if the Coxeter diagram of $(W,S)$ has no connected component of affine
type.
\end{corintro}

The proof of \fullref{thm:FiniteConjug:2sph} makes an essential use of
the $\CAT 0$ cube complex $\mathscr{X}$ constructed by G~Niblo and
L~Reeves \cite{NR03}. More precisely, one associates a cube of
$\mathscr{X}$ to every $2$--spherical subgroup of $W$ in such a way
that the problem of counting conjugacy classes of $2$--spherical
subgroups in $W$ becomes a matter of counting orbits of cubes in
$\mathscr{X}$. In particular, if $W$ acts co-compactly on
$\mathscr{X}$, then there are finitely many orbits of cubes and this
implies immediately that $W$ possesses finitely many conjugacy classes
of $2$--spherical subgroups. However, it is known that the $W$--action
on $\mathscr{X}$ is not always co-compact: as shown in Caprace and
M{\"u}hlherr \cite{CM05}, it is co-compact if and only if $W$ has no
parabolic subgroup of irreducible affine type and rank~$\geq 3$. In
fact, \fullref{thm:FiniteConjug:2sph}(i) can be viewed as a
significant generalization of \cite[Theorem~1.1]{CM05}, but the
somewhat lengthy case-by-case discussions of \cite{CM05} are here
completely avoided. The problem of counting orbits of cubes in the
general case (ie when $\mathscr{X}/W$ is not compact) is settled here
by \fullref{thm:cubes} below, whose proof leads naturally to consider
nested sequences of half-spaces of $\mathscr{A}$. Such sequences are
submitted to a strong alternative described in \fullref{thm:chains}
below.

\subsection{Separation of parallel walls}

A well-known result on Coxeter groups is the so-called \emph{parallel wall theorem}: it asserts that any point
of the Davis complex, which is sufficiently far apart from a given wall, is separated from this wall by another
wall. This was first proved by B~Brink and R~Howlett \cite[Theorem~2.8]{BH94}.
Here, we obtain the following:

\begin{thmintro}\label{thm:ParallelWalls}
There exists a constant $Q=Q(W,S)$, depending only on the Coxeter system $(W,S)$, such that the following holds.
Given two walls $\mu_1, \mu_2$ of $\mathscr{A}$ such that the distance from $\mu_1$ to $\mu_2$ is at least $Q$,
there exists a wall $m$ which separates $\mu_1$ from $\mu_2$.
\end{thmintro}

This was stated in \cite{NR03} as the \emph{strong parallel wall conjecture}. As noted by G~Niblo and
L~Reeves, it implies the existence of a universal bound on the size of the link of a vertex in the $\CAT0$ cube
complex $\mathscr{X}$ constructed in [loc.~cit.]. In fact, the latter corollary can be viewed as an immediate
consequence of the more general fact, due to F~Haglund and D~Wise \cite{HW06}, that the cube complex of any
finitely generated Coxeter group embeds virtually equivariantly in the Davis complex of a right-angled Coxeter
group; the proof of this result relies in an essential way on \fullref{thm:ParallelWalls}.

The proof of this theorem can be outlined as follows. Any two points of $\mathscr{A}$ which are far apart, are
separated by a large set of pairwise parallel walls. In particular, two parallel walls which are far apart yield
a large set of pairwise parallel walls. Each of these walls is a separator-candidate. If it turns out that none
of these candidates indeed separates the original pair of walls, then one constructs a configuration of nested
triangles of walls. The existence of such a configuration is severely restricted by \fullref{thm:ladder}
below.

\subsection{Chains of roots}

As mentioned in the preceding subsections, a common feature of the proofs of \fullref{thm:FiniteConjug:2sph}
and \fullref{thm:ParallelWalls} is that they both lead one to consider at some point sequences of pairwise
parallel walls. In order to study these sequences, we use the  standard root system $\Phi = \Phi(W,S)$. Recall
that the set of roots $\Phi$ is a discrete subset of some finite-dimensional real vector space $V$, which is
endowed with a $W$--action leaving $\Phi$ invariant. This action also preserves a bilinear form $(\cdot, \cdot)$,
called \textit{Tits bilinear form} (see Bourbaki \cite{Bo81} and \fullref{sect:RootBas} below).

This root system is a useful supplement of the Davis complex $\mathscr{A}$, which is a complete $\CAT0$--space
\cite{Da98}, on which $W$ acts properly discontinuously and co-compactly. This complex is a thickened version of
the Cayley graph of $(W,S)$; in fact this Cayley graph is nothing but the $1$--skeleton of $\mathscr{A}$. Walls
of $\mathscr{A}$ are fixed point sets of reflections of $W$. Every wall cuts $\mathscr{A}$ into two open convex
subsets, the closure of which are called half-spaces. The set of all half-spaces (resp.\ walls) is denoted by
$\Phi(\mathscr{A})$ (resp.\ $\mathscr{M(A)}$).

Although different in nature, the root system and the Davis complex are closely related as follows. Once the
identity $1 \in W$ has been chosen as a base vertex of $\mathscr{A}$, one has a canonical $W$--equivariant
bijection $\Phi \to \Phi(\mathscr{A})$, which maps $\Pi$ to the set of those half-spaces containing the vertex
$1$ but not its neighbors (see \fullref{lem:LinkRootBases} below for more details). For the rest of this
introduction, we identify $\Phi$ and $\Phi(\mathscr{A})$ by means of this bijection; thus the words `root' and
`half-space' become synonyms.

\begin{thmintro}\label{thm:chains}
There exists a non-decreasing sequence $(r_n)_{n \in \N}$ of positive real numbers, tending to $+\infty$ with
$n$, such that $r_1 > 1$, which depends only on the Coxeter system $(W, S)$ and such that the following property
is satisfied. Given any chain of half-spaces $\alpha_0 \subsetneq \alpha_1 \subsetneq \dots \subsetneq \alpha_n$
with $n>0$, exactly one of the following alternatives holds:
\begin{itemize}
\item[(1)] $(\alpha_0, \alpha_n) \geq r_n$.

\item[(2)] $(\alpha_0, \alpha_n) = 1$, the group $\la r_{\alpha_i} \; | \; i = 0, \dots, n \ra$ is infinite
dihedral and it is contained in a parabolic subgroup of irreducible affine type of $W$.
\end{itemize}
\end{thmintro}

\fullref{thm:chains} may be viewed as summing up the technical heart of this paper. We will first establish
a weak version of this theorem (\fullref{prop:FiniteChain}), which is sufficient to obtain
\fullref{thm:ParallelWalls}. The latter result is then used to deduce \fullref{thm:chains} in its full
strength (\fullref{sect:chainsII}). This in turn is an essential tool in the proof of
\fullref{thm:FiniteConjug:2sph}.

\section{Preliminaries}

\subsection{Root bases}\label{sect:RootBas}
A \textit{root basis} is a triple $\mathscr{E} = (V,  ( \cdot, \cdot ), \Pi)$ consisting of a real vector space
$V$, a symmetric bilinear form $( \cdot, \cdot )$ on $V$ and a set $\Pi \subset V$ which satisfies the following
conditions:
\begin{itemize}\leftskip 10pt
\item[\textbf{(RB1)}] For all $\alpha \in \Pi$, one has $(\alpha, \alpha)= 1$.

\item[\textbf{(RB2)}] For all $\alpha \neq \beta \in \Pi$, one has
$$(\alpha, \beta) \in \bigl\{-\cos\bigl(\tfrac{\pi}{m}\bigl) \big| \, m \in
\Z_{\geq 2} \bigr\} \cup (-\infty, -1].$$

\item[\textbf{(RB3)}] There exists a linear form $x \in V^*$ such that $x(\alpha) >0$ for all $\alpha \in \Pi$.
\end{itemize}

The most important example is the \textit{standard root basis} of a Coxeter matrix $M=(m_{ij})_{i, j \in I}$. Recall
from \cite{Bo81} that this root basis $\mathscr{E}_M :=(V_M, (\cdot, \cdot)_M, \Pi_M)$ is constructed as
follows: set $V_M:=\bigoplus_{i \in I} \R e_i$, $\Pi_M := \{e_i \; | \; i \in I\}$ and for all $i, j \in I$, set
$(e_i, e_j)_M := -\cos(\frac{\pi}{m_{ij}})$ ($=-1$ if $m_{ij}=\infty$).

Let now $\mathscr{E} = (V,  ( \cdot, \cdot ), \Pi)$ be any root basis. For each $\alpha \in V$, define $r_\alpha
\in \GL(V)$ by
$$r_\alpha : x \mapsto x - 2 (x, \alpha) \alpha.$$
We make the following definitions:
\begin{itemize}
\item $S(\mathscr{E}):= \{r_\alpha | \; \alpha \in \Pi\}$,

\item $W(\mathscr{E}):= \la S(\mathscr{E}) \ra \subset \GL(V)$,

\item $\Phi(\mathscr{E}):= \{w.\alpha | \; w \in W(\mathscr{E}), \; \alpha \in \Pi\}$,

\item $\Phi(\mathscr{E})^+:= \{\phi \in \Phi(\mathscr{E}) | \; \phi \in \sum_{\pi \in \Pi} \R^+ \pi \}$.
\end{itemize}

The elements of $\Phi(\mathscr{E})$ are called \textit{roots}; the roots contained in  $\Phi(\mathscr{E})^+$ are
called \textit{positive}. The following lemma collects the basic facts on root bases which we will need in the
sequel:

\begin{lem}\label{lem:Bourbaki}$\phantom{99}$
\begin{itemize}
\item[(i)] Given $\mathscr{E} = (V,  ( \cdot, \cdot ), \Pi)$ a root basis, the pair $(W(\mathscr{E}),
S(\mathscr{E}))$ is a Coxeter system and $\Phi(\mathscr{E})$ is a discrete subset of $V$.

\item[(ii)] Conversely, if $\mathscr{E}$ is the standard root basis associated with a given a Coxeter system
$(W,S)$, then there is a canonical isomorphism $W \to W(\mathscr{E})$ mapping $S$ onto $S(\mathscr{E})$.

\item[(iii)] For all $w \in W(\mathscr{E})$ and all $\alpha \in \Phi(\mathscr{E})$, one has $w r_\alpha w\inv =
r_{w.\alpha}$.

\item[(iv)] For all $\phi \in  \Phi(\mathscr{E})$, either $\phi \in \Phi(\mathscr{E})^+$ or $-\phi \in
\Phi(\mathscr{E})^+$.
\end{itemize}
\end{lem}
\begin{proof}
See \cite[Theorem~IV.1.1]{Bo81}.  The proofs given in \cite{Bo81} deal
only with standard root bases, but they apply without modification to
any root basis.
\end{proof}
The \textit{type} of the root basis $\mathscr{E}$ is the type of the Coxeter system $(W(\mathscr{E}),
S(\mathscr{E}))$, ie, the Coxeter matrix $M(\mathscr{E}):=(m_{\alpha \beta})_{\alpha, \beta \in \Pi}$, where
$m_{\alpha \beta}$ is the order of $r_\alpha r_\beta \in W(\mathscr{E})$. The following lemma, whose proof is
straightforward, recalls the relationship between two root bases of the same type:

\begin{lem}\label{lem:RootBasis:2sph}
Let $M=(m_{ij})_{i, j \in I}$ be a Coxeter matrix, $\mathscr{E}_M =(V_M, (\cdot, \cdot)_M, \{e_i | \; i \in
I\})$ be the standard root basis of type $M$ and $\mathscr{E} = (V,  ( \cdot, \cdot ), \{\alpha_i | \; i \in
I\})$ be any root basis of type $M$ such that $(\alpha_i, \alpha_j)=-\cos(\pi/m_{ij})$ for all $i, j \in I$ with
$m_{ij} < \infty$. Let $\varphi : W(\mathscr{E}_M) \to W(\mathscr{E})$ be the unique isomorphism such that
$\varphi : r_{e_i} \mapsto r_{\alpha_i}$ for all $i \in I$. We have the following:
\begin{itemize}
\item[(i)] There is a unique $\varphi$--equivariant bijection $f \co  \Phi(\mathscr{E}_M) \to \Phi(\mathscr{E})$
such that $f \co  e_i \mapsto \alpha_i$ for all $i \in I$.

\item[(ii)] If $M$ is $2$--spherical (ie $m_{ij} < \infty$ for all $i, j \in I$), then one has $(\alpha,
\beta)_M = (f(\alpha), f(\beta))$ for all $\alpha, \beta \in \Phi(\mathscr{E}_M)$.\qed
\end{itemize}
\end{lem}

\subsection{Convention}

\emph{ From now on and until the end of the paper, we fix a root basis $\mathscr{E} = (V,  ( \cdot, \cdot ),
\Pi)$ and we set $(W,S) :=(W(\mathscr{E}), S(\mathscr{E}))$. Moreover, given $s \in S$, we denote the unique
element $\alpha$ of $\Pi$ such that $s = r_\alpha$ by $\alpha_s$. }

\subsection{The Davis complex}\label{sect:Davis}

Suppose that the set $S$ is finite. The \textit{Davis complex} $\mathscr{A}$ associated with $(W,S)$ is a
piecewise Euclidean $\CAT0$ cell complex whose $1$--skeleton is the Cayley graph of $W$ with respect to the
generating set $S$. The action of $W$ on this Cayley graph induces naturally an action on $\mathscr{A}$; this
action is properly discontinuous and cocompact. By definition, a \textit{wall} of $\mathscr{A}$ is the fixed
point set of a reflection of $W$. Hence a wall is a closed convex subset of $\mathscr{A}$. A fundamental fact is
that every wall cuts $\mathscr{A}$ into two convex open subsets, whose respective closures are called
\textit{half-spaces}. Hence the boundary $\partial h$ of a half-space $h$ is a wall. The set of all half-spaces
is denoted by $\Phi(\mathscr{A})$.

Given a point $x \in \mathscr{A}$ which does not lie on any wall (e.g. $x$ is a vertex of $\mathscr{A}$), the
intersection $C(x) \subset \mathscr{A}$ of all half-spaces containing $x$ is compact. The set $C(x)$ is called a
\textit{chamber}. Every chamber contains exactly one vertex of $\mathscr{A}$. Hence the $W$--action on
$\mathscr{A}$ is simply transitive on the chambers. Given two chambers $C_1, C_2$, we define the
\textit{numerical distance} from $C_1$ to $C_2$ as the number of walls which separate $C_1$ from $C_2$. If $w
\in W$ is the unique element such that $w.C_1 = C_2$, then this distance equals the word length $\ell(w)$ of
$w$.

Since the $1$--skeleton of $\mathscr{A}$ is the Cayley graph of $(W,S)$, the edges of $\mathscr{A}$ are labelled
by the elements of $S$. Given $s \in S$, two chambers are called \textit{$s$--adjacent} if they contain vertices
which are joined by an edge labelled by $s$.

\begin{lem}\label{lem:LinkRootBases}
Let $C \subset \mathscr{A}$ be a chamber. For each $s \in S$, let $h_s$ be the half-spaces containing $C$ but
not the chamber of $\mathscr{A}$ different from $C$ and $s$--adjacent to $C$. Then one has the following:
\begin{itemize}
\item[(i)] There exists a unique $W$--equivariant bijection $\zeta_C \co  \Phi(\mathscr{E}) \to \Phi(\mathscr{A})$
which maps $\alpha_s$ to $h_s$ for all $s \in S$. The positive roots are mapped onto the half-spaces containing
$C$.

\item[(ii)] Let $\alpha, \beta \in \Phi(\mathscr{E})$. If $|(\alpha, \beta)| < 1$ then the walls $\partial
\zeta_C(\alpha)$ and $\partial \zeta_C(\beta)$ meet. Conversely, if the walls $\partial \zeta_C(\alpha)$ and
$\partial \zeta_C(\beta)$ meet, then $|(\alpha, \beta)| \leq 1$ and equality occurs if and only if $\alpha = \pm
\beta$.

\item[(iii)] For all $\alpha, \beta \in \Phi(\mathscr{E})$, one has $(\alpha, \beta) \geq 1$ if and only if
$\zeta_C(\alpha) \subset \zeta_C(\beta)$ or $\zeta_C(\beta) \subset \zeta_C(\alpha)$.

\item[(iv)] For all $\alpha, \beta \in \Phi(\mathscr{E})$, one has $(\alpha, \beta) \leq 0$ if and only if
$\zeta_C(\alpha) \cap \zeta_C(\beta) \subset \zeta_C(r_\alpha(\beta)) \cap \zeta_C(r_\beta(\alpha))$.
\end{itemize}
\end{lem}
\begin{proof}
Assertion~(i) follows from the fact that the Cayley graph of $(W,S)$ (and even the whole Davis complex) can be
embedded in the Tits cone of the root basis $\mathscr{E}$, see \cite[Appendix~B.4]{Kr94} for details. For (ii)
and (iii), see \cite[Proposition~1.4.7]{Kr94}. Assertion~(iv) follows also by considering the Tits cone.
\end{proof}

\subsection{The cube complex of G~Niblo and L~Reeves}\label{sect:NR}
Maintain the assumption that $S$ is finite. In \cite{NR03}, G~Niblo and L~Reeves used the structure of wall
space of the Cayley graph (or the Davis complex) of $(W,S)$ to construct a $\CAT 0$ cube complex $\mathscr{X}$
endowed with a properly discontinuous $W$--action. We briefly recall here the construction and  basic properties
of $\mathscr{X}$ for later reference.

Vertices of this cube complex are mappings $v\co  \mathscr{M(A)} \to \Phi(\mathscr{A})$ which satisfy the following
two conditions:
\begin{itemize}
\item For all $m \in \mathscr{M(A)}$, we have $\partial(v(m)) = m$.

\item For all $m, m' \in \mathscr{M(A)}$, if $m$ and $m'$ are parallel then either $v(m) \subset v(m')$ or
$v(m') \subset v(m)$.
\end{itemize}
By definition, two distinct vertices are adjacent if and only if the subset of $\mathscr{M(A)}$ on which they
differ is a singleton. Note that the Cayley graph of $(W,S)$ is a subgraph of the so-obtained graph. By
definition, the $1$--skeleton of $\mathscr{X}$ is the connected component of this graph which contains the Cayley
graph of $(W,S)$. The cubes are defined by `filling in' all the cubical subgraphs of the $1$--skeleton. It is
shown in \cite{NR03} that the cube complex $\mathscr{X}$ is finite-dimensional, locally finite, and the
canonical $W$--action is properly discontinuous, but not always co-compact. In fact \cite{CM05}, the $W$--action
on $\mathscr{X}$ is co-compact if and only if $W$ has no parabolic subgroup of irreducible affine type and
rank~$\geq 3$. We refer to \fullref{sect:FiniteConjug} for more precise information on the $W$--orbits of cubes
in $\mathscr{X}$.

As it is the case for any $\CAT 0$ cube complex (see \cite{Sageev}), the space $\mathscr{X}$ is endowed with a
collection $\mathscr{M(X)}$ of \emph{walls} (resp.\ a collection $\Phi(\mathscr{X})$ of \emph{half-spaces}),
which is by construction in canonical one-to-one correspondence with $\mathscr{M(A)}$ (resp.\
$\Phi(\mathscr{A})$). More precisely, a wall is an equivalence class of edges, for the equivalence relation
defined as the transitive closure of the relation of being opposite edges in some square. Thus every edge
defines a wall which separates its two extremities. Given an equivalence class of edges, the corresponding wall
can be realized geometrically as the convex closure of the set of midpoints of edges in this class. In this way,
every wall becomes a closed convex subset of $\mathscr{X}$ which separates $\mathscr{X}$ into two convex
subsets, called half-spaces. Note that a wall is itself a $\CAT 0$ cube complex, which is a subcomplex of the
first barycentric subdivision of $\mathscr{X}$.

\subsection{Reflection subgroups}\label{sect:ReflSubgps}

A \textit{reflection subgroup} of $W$ is a subgroup  generated by some set of reflections. The following basic
fact is well-known:

\begin{lem}\label{lem:ReflSubgps}
Let $H$ be a subgroup of $W$ generated by a set $R$ of reflections. Then there exists a unique set $\Pi' \subset
\Phi(\mathscr{E})^+$ such that $\mathscr{E}'=(V, (\cdot, \cdot), \Pi')$ is a root basis and $W(\mathscr{E}') =
H$. Moreover one has $|\Pi'| \leq |R|$.
\end{lem}
\begin{proof}
See \cite{Deo89} or \cite{Dyer90}.
\end{proof}

By definition, the \textit{type} (resp.\ \textit{rank}) of the reflection group $H$ is the type (resp.\ rank) of
the Coxeter system $(H, S(\mathscr{E}'))$ (see \fullref{lem:Bourbaki}(i)). The reflection group $H$ is called
\textit{standard parabolic} if $S(\mathscr{E}') \subset S$. Let $S(\mathscr{E}') = S_1 \cup \dots \cup S_k$ be
the finest partition of $S(\mathscr{E}')$ into non-empty mutually centralizing subsets. The subgroups $\la S_i
\ra \subset H$, $i=1, \dots, k$, are called the \textit{direct components} of $H$.

\section{Chains of roots}

Throughout this section, we fix a base chamber $C \subset \mathscr{A}$.

\subsection{A partial ordering on the set of roots}\label{sect:kappa}
Transforming the relation of inclusion $\subset$ on $\Phi(\mathscr{A})$ by the bijection $\zeta_C$ of
\fullref{lem:LinkRootBases}(i), one obtains a partial ordering on $\Phi(\mathscr{E})$ which we also denote by
$\subset$. By \fullref{lem:LinkRootBases}, two roots $\alpha, \beta$ are orderable by $\subset$ if and only if
$(\alpha, \beta ) \geq 1$.

Before stating the next lemma, we need to introduce a constant $\kappa$ which is defined as follows:
$$ \kappa = \sup \{ |(\alpha, \beta)| \; : \; \alpha, \beta \in \Phi(\mathscr{E}),
|(\alpha, \beta)| < 1 \}.$$ By \fullref{lem:LinkRootBases}(ii), the condition $|(\alpha, \beta)| < 1$ implies
that the group $\la r_\alpha, r_\beta \ra$ is finite. Since $W$ has finitely many conjugacy classes of finite
subgroups, it follows in particular that $\kappa < 1$. Important to us will be the following:

\begin{lem}\label{lem:OrderedTriples}
Let $\alpha, \beta, \gamma \in \Phi(\mathscr{E})$ be roots such that $\alpha \subset \beta \subset \gamma$. Then
the following holds:
\begin{itemize}
\item[(i)] One has $$(\alpha, \gamma) \geq \max \{ (\alpha, \beta), (\beta, \gamma)\}.$$

\item[(ii)] If moreover $(r_\beta(\alpha), \gamma) > -1$, then $$(\alpha, \gamma) \geq 2\,(\alpha, \beta) -
\kappa.$$

\item[(iii)] If $(r_\beta(\alpha), \gamma) \leq -1$, then $\beta \subset -r_\beta(\alpha) \subset \gamma$ or
$\alpha \subset -r_\beta(\gamma) \subset \beta$.
\end{itemize}
\end{lem}
\begin{proof}
For (i), see \cite[Corollary~6.2.3]{Kr94}. Since $\beta \subset \gamma$,  \fullref{lem:LinkRootBases}(iii)
yields $(\beta, \gamma) \geq 1$. Therefore, one has:
$$\begin{array}{rcl}
(\alpha, \gamma) - 2 \, (\alpha, \beta) &\geq& (\alpha, \gamma) - 2 \, (\alpha, \beta) (\beta, \gamma)\\
& = & (\alpha - 2 \, (\alpha, \beta) \beta,  \gamma)\\
& = & (r_\beta(\alpha), \gamma).
\end{array}$$
Moreover, if $(r_\beta(\alpha), \gamma) > -1$, then $(r_\beta(\alpha), \gamma) \geq - \kappa$ by definition.
This implies (ii). Assertion~(iii) is a consequence of \fullref{lem:LinkRootBases}(ii) and (iii).
\end{proof}

\subsection{An `affine versus non-affine' alternative for chains of roots}
In this section, we establish \fullref{prop:FiniteChain}, which is a first approximation of
\fullref{thm:chains}.

Let $\alpha_1, \alpha_2 \in \Phi(\mathscr{E})$ be roots and let $h_i = \zeta_C(\alpha_i)$ for $i=1,2$. Suppose
that $\alpha_1 \subset \alpha_2$. The set
$$\Phi(\alpha_1 ; \alpha_2):=\{\beta \in \Phi(\mathscr{E}) \; | \; \alpha_1 \subset \beta \subset \alpha_2\}$$
is finite. Indeed, its cardinality is bounded by the combinatorial distance from a vertex contained in $h_2$ to
a vertex contained in the complement of $h_1$.

A set of roots $\Phi \subset \Phi(\mathscr{E})$ is called \textit{convex} if for all $\alpha_1, \alpha_2 \in
\Phi$ and all $\beta \in \Phi(\mathscr{E})$, one has  $\beta \in \Phi$ whenever $\alpha_1 \subset \beta \subset
\alpha_2$. A set of roots $\Phi \subset \Phi(\mathscr{E})$ is called a \textit{chain} if it is totally ordered
by $\subset$. 
A chain is called \emph{convex} if any two consecutive elements form
a convex pair.
In view of the preceding paragraph, it is easy to see that \emph{any chain is contained in a
convex chain}. This convex chain need not be unique, and it is in general properly contained in the convex
closure of the initial chain. A chain of roots $\alpha_0 \subsetneq \alpha_1 \subsetneq \dots \subsetneq
\alpha_n$ is called \textit{maximally convex} if for all chain $\beta_0 \subsetneq \beta_1 \subsetneq \dots
\subsetneq \beta_k$ such that $\beta_0 = \alpha_0$ and $\beta_k = \alpha_n$, one has $k \leq n$. Note that a
maximally convex chain is convex.

As a first consequence of \fullref{lem:OrderedTriples}, we have:
\begin{lem}\label{lem:Phi(alpha;beta)}
Let $\alpha, \beta \in \Phi(\mathscr{E})$ be such that $\alpha \subsetneq \beta$. If $(\alpha, \beta) = 1$, then
the group $\la r_\phi \; | \; \phi \in \Phi(\alpha; \beta) \ra$ is infinite dihedral; in particular
$\Phi(\alpha; \beta)$ is  a chain.
\end{lem}
\begin{proof}
Recall from \fullref{lem:LinkRootBases}(iii) that for all $\phi, \psi \in \Phi(\mathscr{E})$, if $\phi \subset
\psi$ then $(\phi, \psi) \geq 1$.

Choose a chain $\alpha = \alpha_0 \subsetneq \alpha_1 \subsetneq \dots \subsetneq \alpha_k = \beta$ of maximal
possible length; this is possible since $\Phi(\alpha; \beta)$ is finite. Moreover, the maximality of the chain
$(\alpha_i)$ implies that this chain is convex. By \fullref{lem:OrderedTriples}(i) we have $(\alpha_i,
\alpha_j)=1$ for all $i, j= 0, 1, \dots, k$. By \fullref{lem:OrderedTriples}(ii), we have
$(r_{\alpha_{i+1}}(\alpha_i), \alpha_{i+2}) \leq -1$ since $(\alpha_i, \alpha_{i+2}) =1$ for each $i$. Hence, by
\fullref{lem:OrderedTriples}(iii), one has $\alpha_{i+1} \subset -r_{\alpha_{i+1}}(\alpha_i) \subset
\alpha_{i+2}$ or $\alpha_{i} \subset -r_{\alpha_{i+1}}(\alpha_{i+2}) \subset \alpha_{i+1}$. By the convexity of
the chain $(\alpha_i)$, we deduce in both cases that $r_{\alpha_{i+1}}(\alpha_i) = -\alpha_{i+2}$ since
$\alpha_i \neq \alpha_{i+1} \neq \alpha_{i+2}$. Since $i$ was arbitrary, it follows from
\fullref{lem:Bourbaki}(iii) that the group $\la r_{\alpha_i} | \; i=0, 1, \dots, k \ra$ is actually generated
by the pair $\{r_{\alpha_0}, r_{\alpha_1}\}$; in particular this group is infinite dihedral.

Let now $\phi$ be any element of $\Phi(\alpha; \beta)$. We have $-r_\alpha(\beta) \subsetneq -r_\alpha(\alpha_1)
\subsetneq \alpha_0 =\alpha \subset \phi \subset \beta$. Since $(\alpha, \beta) = 1$, it follows that
$(-r_\alpha(\beta), \beta)=1$. Hence, by \fullref{lem:OrderedTriples}(i) we have $(-r_\alpha(\alpha_1), \phi)
= 1$. Therefore, the restriction of the bilinear form $(\cdot, \cdot)$ to the subspace spanned by
$\{-r_\alpha(\alpha_1), \alpha, \phi\}$ is positive semi-definite and that its radical is of codimension~$1$. By
\fullref{lem:ReflSubgps} and \cite[Chapter~VI, Section~4.3, Theorem~4]{Bo81}, this implies that the reflection subgroup
generated by $\{r_\alpha, r_{\alpha_1}, r_\phi\}$ is infinite dihedral. By the above, the latter group contains
$r_{\alpha_i}$ for each $i$. Therefore the wall $\partial \zeta_C(\phi)$ is parallel to $\partial
\zeta_C(\alpha_i)$ for each $i$. By the maximality of the chain $(\alpha_i)_{i \leq k}$, this implies that $\phi
= \alpha_j$ for some $j \in \{0, 1, \dots, k\}$. In other words, we have $\Phi(\alpha; \beta) = \{\alpha_0,
\alpha_1, \dots, \alpha_k\}$, which completes the proof.
\end{proof}

A posteriori, the preceding lemma can be viewed as a consequence of \fullref{thm:chains}; however, the proof
of the latter relies on \fullref{lem:Phi(alpha;beta)} in an essential way.

The constant $\kappa$ appearing in the next proposition was defined in \fullref{sect:kappa}.

\begin{prop}\label{prop:FiniteChain} Let $\alpha_0 \subsetneq \alpha_1 \subsetneq \dots \subsetneq
\alpha_n$ be a maximally convex chain of roots and let $j \in \{1, \dots, n\}$. We have the following:
\begin{itemize}
\item[(i)] Assume that $(\alpha_0, \alpha_{j-1}) = 1$. Then either $(\alpha_0, \alpha_j) = 1$ or $(\alpha_0,
\alpha_j) \geq j (1 - \kappa)$.

\item[(ii)] Assume that $(\alpha_0, \alpha_j) = 1 + \varepsilon$ for some $\varepsilon > 0$. Then $(\alpha_0,
\alpha_n) > 1 + \frac{n}{2j}\varepsilon $.
\end{itemize}
\end{prop}
\begin{proof}
(i)\qua If $j= 1$, there is nothing to prove. Thus we assume $j > 1$. By \fullref{lem:Phi(alpha;beta)} and by
convexity of the chain $(\alpha_i)$, the group $\la r_{\alpha_i} | \; i \in \{0, 1, \dots, j-1\}\ra$ is infinite
dihedral, generated by the pair $\{r_{\alpha_{j-2}}, r_{\alpha_{j-1}}\}$. Let $\beta =
r_{\alpha_{j-1}}(\alpha_{j-2})$.

Assume first that $(\beta, \alpha_j) \leq -1$. Then \fullref{lem:OrderedTriples}(iii)  implies $\beta = -
\alpha_j$ and, hence, $r_{\alpha_j} \in \la r_{\alpha_i} | \; i \in \{0, 1, \dots, j-1\}\ra$. An easy
computation using \fullref{lem:Bourbaki}(iii) yields $(\alpha_0, \alpha_j) = 1$.

Assume now that $(\beta, \alpha_j) > -1$. Clearly this implies $(\beta, \alpha_j) \geq -\kappa$. Using
\fullref{lem:Bourbaki}(iii), one easily computes that $(\beta, \alpha_{j-1}) = -1$ and that $\alpha_0 = j.
\alpha_{j-1} + (j-1) . \beta$. Therefore, we deduce
$$\begin{array}{rcl}
(\alpha_0, \alpha_j) &=& j.(\alpha_{j-1}, \alpha_j) + (j-1).(\beta, \alpha_j)\\
&\geq& j - (j-1)\kappa \\
&=& j (1 -\kappa) + \kappa.
\end{array}$$

(ii)\qua By \fullref{lem:OrderedTriples}(i), we have $(\alpha_0, \alpha_n) \geq (\alpha_0, \alpha_j)$.
Thus we may assume $n \geq 2j$, otherwise there is nothing to prove. Clearly $\alpha_j \subsetneq
-r_{\alpha_j}(\alpha_0)$. Let
$$k = \max \{ i \; | \; \alpha_i
\subsetneq -r_{\alpha_j}(\alpha_0) \}.$$ Thus $k \geq j$.

Assume that $k \geq 2j$, namely that $\alpha_{2j} \subsetneq -r_{\alpha_j}(\alpha_0)$. Then one would have a
chain
$$\alpha_0 \subsetneq -r_{\alpha_j}(\alpha_{2j}) \subsetneq -r_{\alpha_j}(\alpha_{2j-1}) \subsetneq \dots
\subsetneq -r_{\alpha_j}(\alpha_{j+1}) \subsetneq \alpha_j$$ contradicting the fact that the sequence
$(\alpha_i)_{i \leq n}$, and hence $(\alpha_i)_{i \leq j}$, is maximally convex. Thus $k < 2j$. Let us now
consider the ordered triple $\alpha_0 \subsetneq \alpha_j \subsetneq \alpha_{k+1}$.

Suppose first that $(r_{\alpha_j}(\alpha_0), \alpha_{k+1}) \leq -1$. Then \fullref{lem:OrderedTriples}(iii)
implies that $-r_{\alpha_j}(\alpha_0) \subset \alpha_{k+1}$ or $\alpha_{k+1} \subset -r_{\alpha_j}(\alpha_0)$.
On the other hand,  the definition of $k$ implies that $\alpha_k \subsetneq -r_{\alpha_j}(\alpha_0)$ and that
$\alpha_{k+1}$ is not properly contained in $-r_{\alpha_j}(\alpha_0)$. Therefore, by the convexity of the chain
$(\alpha_i)_{i \leq n}$, we have $r_{\alpha_j}(\alpha_0) = -\alpha_{k+1}$ in this case. We deduce that
$$\begin{array}{rcl}
(\alpha_0, \alpha_{k+1}) &=& (\alpha_0, -r_{\alpha_j}(\alpha_0))\\
& = &  2(\alpha_0, \alpha_j)^2 -1\\
& = &  1 + 4 \varepsilon + 2 \varepsilon^2\\
& > &  1 + 2 \varepsilon.
\end{array}$$
Suppose now that $(r_{\alpha_j}(\alpha_0), \alpha_{k+1}) > -1$. Then \fullref{lem:OrderedTriples}(ii) implies
that $(\alpha_0, \alpha_{k+1}) \geq 2 (\alpha_0, \alpha_j)-\kappa > 1 + 2 \varepsilon$.

In both cases, we obtain $(\alpha_0, \alpha_{2j}) > 1 + 2 \varepsilon$ by \fullref{lem:OrderedTriples}(i)
because $k+1 \leq 2j$. An immediate induction now yields $(\alpha_0, \alpha_{2^x j}) > 1 + 2^x \varepsilon$ for
all positive integer $x$ such that $2^x j\leq n$. Since the maximal such integer is $\lfloor \log_2(\frac{n}{j})
\rfloor$, we deduce, again from \fullref{lem:OrderedTriples}(i), that $(\alpha_0, \alpha_n) > 1+2^x
\varepsilon$ with $x =\lfloor \log_2(\frac{n}{j}) \rfloor$. The desired inequality follows because $\lfloor
\log_2(\frac{n}{j}) \rfloor > \log_2(\frac{n}{j})-1 = \log_2(\frac{n}{2j})$.
\end{proof}

\section{Nested triangles of walls}

When studying Coxeter groups, it is often useful to relate the combinatorics of walls in the Davis complex with
the algebraic properties of the subgroup generated by the corresponding reflections. A typical well-known result
of this kind is the basic fact that two walls meet if and only if the corresponding reflections generate a
finite group. The purpose of this section is prove the following:

\begin{thm}\label{thm:ladder}
There exists a constant $L=L(W,S)$, depending only on the Coxeter system $(W,S)$, such that the following
property holds. Let $\mu, \mu'$, $m_0, m_1, \dots, m_n$ be walls of the Davis complex $\mathscr{A}$ such that:
\begin{itemize}
\item[(1)] $\varnothing \neq \mu \cap \mu' \subset m_0$;

\item[(2)] For all $0 \leq i < j < k \leq n$, the wall $m_j$ separates $m_i$ from $m_k$;

\item[(3)] For each $i = 1, \dots, n$, the wall $m_i$ meets both $\mu$ and $\mu'$.
\end{itemize}
If $n > L$, then the group generated by the reflections through the walls $\mu$, $\mu'$, $m_0$, $m_1, \dots,
m_n$ is isomorphic to a Euclidean triangle group. Moreover, it is contained in a parabolic subgroup of
irreducible affine type.
\end{thm}

By a \textit{Euclidean triangle group}, we mean an affine Coxeter group of rank~$3$, or equivalently, the
automorphism group of one of the three (types of) regular tessellations of the Euclidean plane by triangles.

\fullref{thm:ladder} has also proved essential in studying Euclidean flats isometrically embedded in the
Davis complex \cite{CH06}. In this paper, it is an essential ingredient in the proof of
\fullref{thm:ParallelWalls}.

\subsection{Nested Euclidean triangles}

Given a set of walls $M$, we denote by $W(M)$ the reflection subgroup of $W$ generated by all reflections
associated to elements of $M$.

\begin{lem}\label{lem:nestedEucl}
Let $\mu, \mu', m_0, m_1, \dots, m_n$ be walls of the Davis complex $\mathscr{A}$ which satisfy conditions (1),
(2) and (3) of \fullref{thm:ladder}. Assume that the group $W(\{\mu, \mu', m_0, m_n\})$ is a Euclidean
triangle subgroup. Then we have the following:
\begin{itemize}
\item[(i)] The group $W(\{\mu, \mu', m_0, m_1, \dots, m_n\})$ is a Euclidean triangle subgroup, which is
contained in a parabolic subgroup of affine type of $W$; in particular\break $W(\{ m_0, m_1, \dots, m_n\})$ is
infinite dihedral.

\item[(ii)] For any point $x \in \mu' \cap m_n$, there exist $k = \lfloor \frac{n}{2} \rfloor$ pairwise parallel
walls\break $m'_1, m'_2, \dots, m'_k$ which separate $x$ from $\mu$, and such that $$W(\{m'_1, m'_2, \dots, m'_k\})
\subset W(\{\mu, \mu', m_0, m_n\});$$ in particular $W(\{\mu, m'_1, m'_2, \dots, m'_k\})$ is infinite dihedral.
\end{itemize}
\end{lem}
\begin{proof}
(i)\qua By a theorem of D~Krammer which is recalled in \fullref{prop:PcAffine} below, the group $W(\{\mu,
\mu', m_0, m_n\})$ is contained in a parabolic subgroup $W_0$ of irreducible affine type of $W$. The Davis
complex of this parabolic subgroup is contained in $\mathscr{A}$ as a residue $\rho_0$; in other words, there is
a chamber $C_0$ such that the union $\rho_0 = \bigcup_{w \in W_0} w.C_0$ is a closed convex subset of
$\mathscr{A}$, whose stabilizer in $W$ coincides with $W_0$. Since the reflection $r_{m_0}$ belongs to $W(\mu,
\mu')$ by condition~(1), it follows that $r_{m_0}$ and $r_{m_n}$ both stabilize $\rho_0$. In particular, the
walls $m_0$ and $m_n$ both meet $\rho_0$. Since $\rho_0$ is convex, it follows that every wall which separates
$m_0$ from $m_n$ meets $\rho_0$. Since $\rho_0$ is a residue, every reflection associated to a wall which cuts
$\rho_0$ must stabilize $\rho_0$. It follows that $r_{m_i} \in \Stab_W(\rho_0) = W_0$ for all $i = 1, \dots, n$,
as desired. Finally, it is an easy observation that any subgroup of an affine Coxeter group generated by
reflections through pairwise parallel walls is infinite dihedral.

(ii)\qua We have just seen that any wall $m$ which separates $m_0$ from $m_n$ meets $\rho_0$ and, hence, belongs to
$W_0$. Such a wall is parallel to $m_i$ for each $i$ because parallelism of walls in an affine Coxeter group is
an equivalence relation. This shows that we may assume, without loss of generality, that every wall which
separates $m_0$ from $m_n$ is one of the $m_i$'s.

By assumption (1) the reflection $r_\mu$ and $r_{\mu'}$ do not commute. Therefore, by considering each of the
three types of affine triangle groups separately, it is easily seen that for each $i = 1, \dots, k = \lfloor
\frac{n}{2} \rfloor$, there exists a wall $m'_i$ which is parallel to $\mu$ and such that $m_{2i} \cap \mu'
\subset m'_i$. Choose any point $y$ on $\mu \cap \mu'$ and consider a geodesic path joining $x$ to $y$. This
path is completely contained in $\mu'$ by convexity, and meets each $m_i$ by assumption (2). Since $m_{2i} \cap
\mu' \subset m'_i$, it follows that this path crosses $m'_i$ for each $i = 1, \dots, k $. Since $m'_i$ is
parallel to $\mu$, this means precisely that $\mu'_i$ separates $x$ from $\mu$.
\end{proof}

\subsection{Critical bounds for hyperbolic triangles}

Before stating the next lemma, we introduce an additional constant $\lambda_{\mathrm{fin}}$ which is defined as
follows:
$$\lambda_{\mathrm{fin}}
= \sup \{ x \in \R \; | \; \alpha = x.\phi + y.\psi, \;  \alpha, \phi, \psi \in \Phi(\mathscr{E}), \; |(\phi,
\psi)| < 1\}.$$ (By convention, we set $\lambda_{\mathrm{fin}} = 1$ if $|(\phi, \psi)| \geq 1$ for all $\phi,
\psi \in \Phi$.) Note that the conditions $ \alpha = x.\phi + y.\psi$ and $|(\phi, \psi)| < 1$ imply that the
group $\la r_\alpha, r_\phi, r_\psi \ra$ is a finite dihedral group. Since $W$ has finitely many conjugacy
classes of finite subgroups, it follows that the constant $\lambda_{\mathrm{fin}}$ is finite; thus
$\lambda_{\mathrm{fin}}$ is a positive real number.

\begin{lem}\label{lem:HypTr}
Let $\phi, \phi', \alpha_0, \alpha_1 \in \Phi(\mathscr{E})$ be roots such that $\alpha_0 \subsetneq \alpha_1$
and that the walls $\mu =
\partial \phi$, $\mu' = \partial \phi'$, $m_0 = \partial \alpha_0$ and $m_1 = \partial \alpha_1$ satisfy
conditions (1) and (3) of \fullref{thm:ladder}. If  $W(\mu, \mu', m_0, m_1)$ is not a Euclidean triangle
subgroup, then one has
$$1 + \varepsilon \leq (\alpha_0, \alpha_1) \leq 2 \kappa \lambda_{\mathrm{fin}},$$ where $\varepsilon =
\varepsilon(W,S) >0$ is a positive constant which depends only on $(W,S)$ and $\kappa$ is the constant defined
in \fullref{sect:kappa}.
\end{lem}

It is well-known that, in the situation of the preceding lemma, if  $W(\mu, \mu', m_1)$ is a Euclidean triangle
subgroup, then $(\alpha_0, \alpha_1)=1$ (see the first lines of the proof of
\fullref{prop:AffinePairs}(ii) below).

\begin{proof}
Condition~(1) implies that $W(\mu, \mu', m_0)$ is a finite dihedral group. In particular we have $\alpha_0 =
\lambda_0.\phi + \mu_0.\phi'$ for some  $\lambda_0, \mu_0 \in \R$. We have
$$(\alpha_0, \alpha_1) =
\lambda_0.(\phi, \alpha_1) + \mu_0.(\phi', \alpha_1) \leq \kappa(\lambda_0 + \mu_0) \leq 2 \kappa
\lambda_{\mathrm{fin}}.$$
Clearly the group $W(\mu, \mu', m_0, m_1)$ is infinite because the walls $m_0$ and $m_1$ are parallel. Since
$W(\mu, \mu', m_0)$ is a finite dihedral group, it follows from \fullref{lem:ReflSubgps} that the reflection
group $W(\mu, \mu', m_0, m_1)$ is infinite of rank~$3$. Since it is not of affine type by hypothesis, it must be
a hyperbolic triangle group, namely it is isomorphic to the automorphism group of a regular tessellation of the
hyperbolic plane by triangles. Furthermore condition~(3) implies that $W(\mu, \mu', m_0, m_1)$ is compact
hyperbolic, ie, the tiles of the above tessellation are compact triangles. We view this tessellation as a
geometric realization of the Coxeter complex associated to $W(\mu, \mu', m_0, m_1)$.

In this realization, the walls $\partial \alpha_0$ and $\partial \alpha_1$ are realized by parallel geodesic
lines. It is known \cite[Corollary~A.5.8]{BenedettiPetronio} that, if one writes $(\alpha_0, \alpha_1)
=\frac{1}{2}(x + x\inv)$ for some $x \geq 1$, then the hyperbolic distance between these lines in $\Hyp^2$ is
$\log(x)$ (see also \fullref{lem:RootBasis:2sph}(ii)). In view of this formula, we have $(\alpha_0, \alpha_1)
> 1$ because parallel walls do not meet in $\Hyp^2 \cup \partial \Hyp^2$ in a regular tessellation by compact
triangles.

It remains to show that $(\alpha_0, \alpha_1)$ stays bounded away from $1$ when $\phi, \phi', \alpha_0$ and
$\alpha_1$ vary in $\Phi(\mathscr{E})$. First, we note that, by the formula above and the fact that $W(\mu,
\mu', m_0, m_1)$ is transitive on pairs of parallel walls at small distances (see \fullref{lem:OrbitsPairs}
below), the scalar $(\alpha_0, \alpha_1)$ stays bounded away from $1$ when $\phi, \phi', \alpha_0$ and
$\alpha_1$ vary in such a way that the group $W(\mu, \mu', m_0, m_1)$ remains in the same isomorphism class.
Now, the desired result follow because $W(\mu, \mu', m_0, m_1)$ is a $2$--spherical reflection subgroup of
rank~$3$ of $W$, and  $W$ has finitely many types of such reflection subgroups: indeed, the type of a such a
subgroup is a Coxeter matrix of size~$3$, all of whose entries divide some entry of the Coxeter matrix of
$(W,S)$. This concludes the proof.
\end{proof}

\begin{remarknb}
\fullref{prop:AffinePairs}(ii) below, which relies on the preceding lemma through
Theorems~\ref{thm:ladder} and~\ref{thm:ParallelWalls}, may be viewed as a generalization of the first inequality
of \fullref{lem:HypTr}.
\end{remarknb}

\begin{remarknb} It turns out that the positive constant $\varepsilon$ appearing in \fullref{lem:HypTr} can be made
`universal', ie independent of $(W,S)$. This is done by expressing the minimal hyperbolic distance between two
walls of a regular tessellation of $\Hyp^2$ by compact triangles, as a monotonic function of the area of the
fundamental tile $T$. This tile is a triangle whose angles are of the form $(\frac{\pi}{k}, \frac{\pi}{l},
\frac{\pi}{m})$ for some $k, l, m \in \N$. Thus the area of $T$ is $\pi - \frac{\pi}{k} - \frac{\pi}{l} -
\frac{\pi}{m}$, and it is not difficult to compute that this area has a global minimum, which is reached for
$(k, l, m )= (2, 3, 7)$. However, the universality of $\varepsilon$ is not relevant to our purposes.
\end{remarknb}

\subsection[Proof of \ref{thm:ladder}]{Proof of \fullref{thm:ladder}}

We will show that the constant $L$ can be defined by
$$L = \max \bigg\{ 2, \ \frac{2\kappa \lambda_{\mathrm{fin}}}{1-\kappa}, \
\frac{8 \kappa^2 \lambda_{\mathrm{fin}}^2 - 4 \kappa \lambda_{\mathrm{fin}} }{ \varepsilon (1-\kappa)}\bigg\}.$$
Let $\phi, \phi', \alpha_0, \alpha_1, \dots, \alpha_n \in \Phi(\mathscr{E})$ be roots such that $\alpha_0
\subsetneq \alpha_1 \subsetneq \dots \subsetneq \alpha_n$ and that the walls $\mu =
\partial \phi$, $\mu' = \partial \phi'$ and $m_i = \partial \alpha_i$ for $i=0, 1, \dots, n$ satisfy
conditions (1), (2) and (3) of \fullref{thm:ladder}. Suppose moreover that $n > L$.

By \fullref{lem:nestedEucl}, we may assume that $W(\mu, \mu', m_0, m_n)$ is not a Euclidean triangle subgroup,
otherwise we are done. Therefore, \fullref{lem:HypTr} yields $1+\varepsilon \leq (\alpha_0, \alpha_n) \leq 2
\kappa \lambda_{\mathrm{fin}}$.

We now make another estimate of the value of $(\alpha_0, \alpha_n)$. To this end, note first that every wall
which separates $m_0 = \partial \alpha_0$ from $m_n = \partial \alpha_n$ must meet $\mu = \partial \phi$ and
$\mu' = \partial \phi'$ because $\mu$ and $\mu'$ are convex, as are all walls. Thus, after replacing the
sequence $(\alpha_i)_{i \leq n}$ by a maximally convex chain of roots whose extremities are $\alpha_0$ and
$\alpha_n$, one obtains a new set of roots which satisfies again conditions (1), (2) and (3) of the theorem. We
henceforth assume without loss of generality that the sequence $(\alpha_i)_{i \leq n}$ is maximally convex.

Let $j = \min \{i \; | \; (\alpha_0, \alpha_i) > 1 \}$. Thus $(\alpha_0, \alpha_{j-1}) = 1$ and
\fullref{prop:FiniteChain}(i) yields $(\alpha_0, \alpha_j) \geq j(1-\kappa)$. By
\fullref{lem:OrderedTriples}(i), we have $(\alpha_0, \alpha_j) \leq (\alpha_0, \alpha_n)$ and, hence, we
deduce that $j \leq \frac{2\kappa \lambda_{\mathrm{fin}}}{1-\kappa}$.

{We now apply \fullref{prop:FiniteChain}(ii). This yields $(\alpha_0, \alpha_n) > 1+\frac{n}{2j}
\varepsilon$. Since $n> L$ we have
$$\frac{n}{2j} \geq \frac{n(1-\kappa)}{4 \kappa \lambda_{\mathrm{fin}} }
> \frac{2 \kappa \lambda_{\mathrm{fin}} -1 }{\varepsilon}.$$
We deduce that $(\alpha_0, \alpha_n) > 2\kappa \lambda_{\mathrm{fin}}$, a contradiction. This shows that the
group\break $W(\mu, \mu', m_0, m_n)$ is a Euclidean triangle subgroup and, hence, the desired result follows from
\fullref{lem:nestedEucl}.}\qed

\section{Separation of parallel walls}

\subsection{On the walls of an infinite dihedral subgroup}

Note that if a set $M$ of walls is such that  the group $W(M)$ is infinite dihedral, then the elements of $M$
are pairwise parallel.

\begin{lem}\label{lem:Rk3}
Let $M$ be a set of walls such that the group $W(M)$ is infinite dihedral. Let $m$ be a wall which meets at
least $8$~elements of $M$ in $\mathscr{A}$. Then $m$ meets all elements of $M$ and either the reflection $r_m$
through $m$ centralizes $W(M)$ or $W(M \cup\{m\})$ is a Euclidean triangle subgroup.
\end{lem}
\begin{proof}
By \fullref{lem:ReflSubgps}, we may -- and shall -- assume, without loss of generality, that $W = W(M \cup
\{m\})$. Since $W$ is infinite, it is of rank~$\geq 2$. Moreover $W$ is not infinite dihedral otherwise the
walls in $M \cup \{m\}$ would be pairwise parallel, in contradiction with the hypotheses. Thus $W$ is of rank~$>
2$. On the other hand, it is generated by a reflection together with a dihedral reflection group; hence it has a
generating set consisting of $3$ reflections. By \fullref{lem:ReflSubgps}, this shows that $W$ is of rank~$3$.

If $W$ is of reducible type, then the reflection $r_m$ through $m$ must centralize $W(M)$ and, hence, the wall
$m$ meets every element of $M$. Thus we are done in this case, and we assume from now on that $W$ is of
irreducible type. A Coxeter group of irreducible type and rank~$3$ is either a Euclidean triangle subgroup, or a
hyperbolic triangle subgroup, ie, a group isomorphic to the automorphism group of a regular tessellation of the
hyperbolic plane $\Hyp^2$ by triangles.

If $W$ is a Euclidean triangle subgroup, then $m$ meets every element of $M$, since parallelism of walls is an
equivalence relation in Euclidean geometry. Thus we are done in this case as well, and it remains to deal with
the case when $W$ is a hyperbolic triangle subgroup. In particular $W$ is Gromov-hyperbolic. We consider the
regular tessellation of $\Hyp^2$ by triangles whose automorphism group is isomorphic to $W$, and view this
tessellation as a geometric realization of $(W,S)$. In particular the walls are geodesic lines in $\Hyp^2$.

Up to enlarging $M$ if necessary, we may assume that for each wall $\mu$, if $r_\mu \in W(M)$ then $\mu \in M$.
We first prove that $m$ meets only finitely many elements of $M$. Suppose the contrary in order to obtain a
contradiction. Choose a chain of half-spaces $\dots \subsetneq h_{-1} \subsetneq h_0 \subsetneq h_1 \subsetneq
\dots$ such that $M = \{ \partial h_i \; | \; i \in \Z\}$ and that for each $i \in \Z$, the set $h_i \cap m$
contains a geodesic ray. This is possible because $m$ crosses infinitely many $\partial h_i$'s. Note that the
intersection $\bigcap_{i \in \Z} h_i$ is a single point $\xi$ of the visual boundary $\partial \Hyp^2$, which
must be an endpoint of $m$ by the above. Thus $\xi$ is fixed by $r_m$ and by the translation subgroup of $W(M)$.
But the stabilizer of $\xi$ in $\Isom(\Hyp^2)$ is solvable. Hence, given a nontrivial translation $t \in W(M)$,
the group $\la t, r_m t r_m \ra$ is solvable, whence virtually abelian of $\Z$--rank~$1$. This shows that $r_m t
r_m$ fixes both of the two points at infinity fixed by $t$, hence so does $r_m$. Furthermore, the pair
consisting of these two points is stabilized by $W(M)$, which finally shows that $W = W(M \cup \{m\})$
stabilizes a pair of points of $\partial \Hyp^2$. This is absurd because $W$ is not virtually solvable. Thus $m$
meets finitely many elements of $M$.

The next step is to show that $r_m$ commutes with at most one reflection in $W(M)$. Indeed, if $\mu$ and $\mu'$
are two distinct elements of $M$ such that $r_m$ centralizes $\la r_\mu, r_{\mu'} \ra$ then $m$ meets every
translate of $\mu$ under $\la r_\mu, r_{\mu'} \ra$, but there are infinitely many such translates, which all
belong to $M$. This is absurd, because we have just seen that $m$ meets finitely many elements of $M$.

We now consider again a chain of half-spaces $\dots \subsetneq h_{-1} \subsetneq h_0 \subsetneq h_1 \subsetneq
\dots$ such that $M = \{ \partial h_i \; | \; i \in \Z\}$. We assume moreover that the numbering of the $h_i$'s
is such that the set of all those $j$'s such that $m$ meets $\partial h_j$ is $\{-k', \dots, 0, \dots, k\}$ for
some integers $k, k'$ with $k' + 1\geq k \geq k' \geq 0$.

We have seen that the angle between $m$ and $\partial h_j$ equals $\frac{\pi}{2}$ for at most one $j \in \Z$. By
the choice of numbering of the $h_i$'s made above, this means in particular that the only possible $j$ such that
$m$ is perpendicular to $\partial h_j$ is $j=0$. In particular the wall $m_1 = r_{\partial h_1}(m)$ is different
from $m$. Since $m \cap \partial h_j \neq \varnothing $ for all $j \in \{-k', \dots, 0, \dots, k\}$ and since
$r_{\partial h_1}(\partial h_i) = h_{2-i}$ for all $i$, it follows that $m_1 \cap \partial h_j \neq \varnothing$
for all $j \in \{2 - k, \dots, 2+k' \}$. This means that the walls $m, m_1$ and $\partial h_{k''}$ form a
compact geodesic triangle $T$ in $\Hyp^2$, where $k'' = \min\{k, \; 2+k'\}$. Since $k \leq k'+1$, we have $k'' =
k$. Furthermore, the triangle $T$ is cut by $\partial h_j$ for all $j \in \{2, \dots, k-1\}$. Now, the
triangular tessellation of $\Hyp^2$ by chambers induces a Coxeter decomposition of $T$ by triangles, namely a
tessellation of $T$ by triangles such that two triangles sharing an edge are switched by the orthogonal
reflection through that edge. But all Coxeter decompositions of hyperbolic triangles are classified
\cite[Section 5.1]{Felikson}. Using this classification, together with the fact that the angle between $m$ and
$\partial h_j$ is $< \frac{\pi}{2}$ for all $j \in \{1, \dots, k\}$, one deduces easily that $k \leq 3$. In
particular, we have also $k' \leq 3$. This shows that $m$ meets at most $7$~elements of $M$.
\end{proof}

\subsection{On walls which separate a vertex from its projection}

\begin{lem}\label{lem:projection}
Let $\phi \in \Phi(\mathscr{A})$ be a half-space, $x \in \mathscr{A}$ be a vertex not contained in $\phi$ and
$y$ be a vertex contained in $\phi$ and at minimal combinatorial distance from $x$. Let $\psi$ be a half-space
containing $x$ but not $y$, such that $\partial \phi \neq \partial \psi$. Then we have the following:
\begin{itemize}
\item[(i)] $\phi \cap \psi \subset r_\phi(\psi) \neq \psi$.

\item[(ii)] $x, y \not \in r_\phi(\psi)$.

\item[(iii)] If moreover $\partial \phi$ meets $\partial \psi$  then $\partial r_\phi(\psi)$ meets every wall
which separates $x$ from $\partial \psi$ and which meets $\partial \phi$.
\end{itemize}
\end{lem}
\begin{proof}
For any half-space $h$ we denote by $-h$ the other half-space bounded by $\partial h$.

If $\phi \cap \psi = \emptyset$, then the walls $\partial \phi$ and $\partial \psi$ are parallel and
$r_\phi(\psi)$ is properly contained in $\phi$. We deduce that $x \not \in r_\phi(\psi)$. By the minimality
assumption on $y$, it follows that $-\phi$ contains a vertex $y'$ neighboring $y$. Note that $y' = r_\phi(y)$.
Since $\partial \phi$ is the only wall separating $y$ from $y'$ we have $y' \in -\psi$, which yields $y =
r_\phi(y') \in r_\phi(-\psi) = - r_\phi(\psi)$, or equivalently $y \not \in r_\phi(\psi)$.

Thus we may assume that $\phi \cap \psi$ is nonempty.  This implies that the walls $\partial \phi$ and $\partial
\psi$ meet, otherwise $\phi$ would be contained in $\psi$ which is impossible since $y \in \phi \cap (-\psi)$ by
assumption. Therefore $\la r_\phi, r_\psi \ra$ is a finite dihedral group.

For any two half-spaces $\alpha, \beta$ such that $\la r_\alpha, r_\beta \ra$ is a finite dihedral group, it is
straightforward to check the following observations (see also \fullref{lem:LinkRootBases}):
\begin{enumerate}
\item[(1)] $\alpha \cap \beta \subset r_\alpha(\beta) \Leftrightarrow \alpha \cap \beta \subset
r_\beta(\alpha)$.

\item[(2)] $\alpha \cap \beta \subset r_\alpha(\beta)  \Leftrightarrow (-\alpha) \cap (-\beta) \subset
-r_\alpha(\beta)$.

\item[(3)] $\alpha \cap \beta \not \subset r_\alpha(\beta) \Rightarrow (-\alpha) \cap \beta \subset
r_\alpha(\beta)$.
\end{enumerate}

Assume now that $\phi \cap \psi \not \subset r_\phi(\psi)$ in order to obtain a contradiction. By (3), this
yields $(-\phi) \cap \psi \subset r_\phi(\psi)$. Since $x \in (-\phi) \cap \psi$, we deduce from (1) that $x \in
-r_\psi(\phi)$; similarly, since $y \in \phi \cap (-\psi)$, we deduce from (1) and (2) that $y \in
r_\psi(\phi)$.

Consider now a minimal combinatorial path $\gamma \co  x = x_0, x_1, \dots, x_n = y$ joining $x$ to $y$ in the
$1$--skeleton of $\mathscr{A}$ (see \fullref{fig:Projection}). We first show that the reflections $r_\phi$ and
$r_\psi$ do not commute. If they commuted, we would have $r_\psi(\phi) = \phi$ and hence, the vertex $r_\psi(y)$
would belong to $\phi$. On the other hand, the path $\gamma$ crosses $\partial \psi$, namely there exists $i \in
\{1, \dots, n\}$ such that $\psi$ contains $x_{i-1}$ but not $x_i$. Thus the path $x_0, \dots, x_{i-1} =
r_\psi(x_i), r_\psi(x_{i+1}), \dots, r_\psi(y)$, which is of length $n-1$, joins $x$ to a vertex of $\phi$,
which contradicts the minimality assumption on $y$. This shows that $r_\phi$ and $r_\psi$ do not commute and,
hence, that $\partial r_\psi(\phi) \neq \partial \phi$.

\begin{figure}[ht!]\label{fig:Projection}
\begin{center}
\includegraphics[scale=0.48]{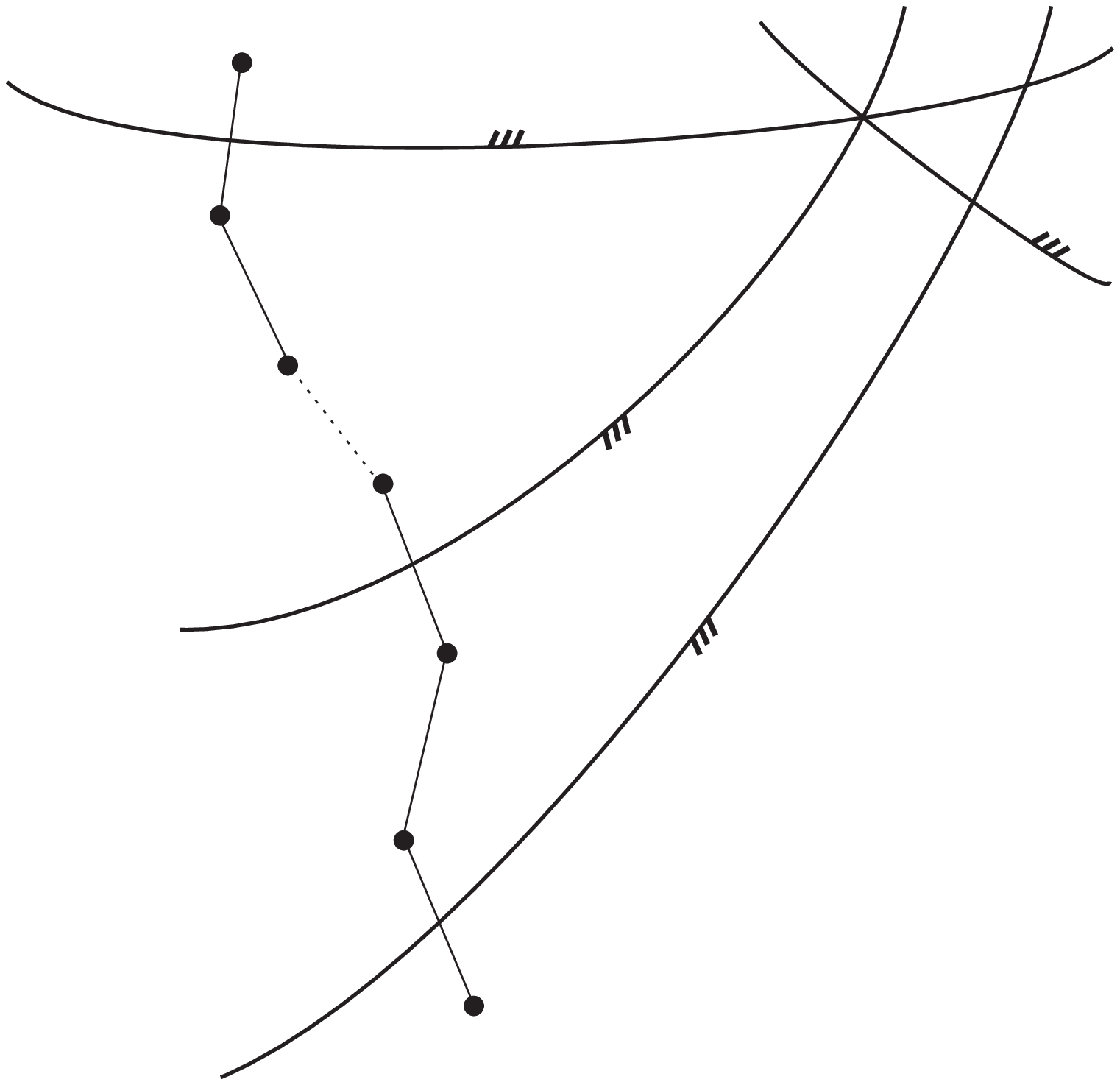}
\unitlength=.8pt\small
\begin{picture}(0,0)
\put(-150,15){$x$}%
\put(-189,55){$x_1$}%
\put(-211,202){$x_{n-1}$}%
\put(-205,237){$y$}%
\put(-150,230){$\phi$}%
\put(-120,141){$\psi$}%
\put(-93,95){$\alpha$}%
\put(-20,204){$r_\phi(\psi)$}%
\end{picture}
\end{center}
\caption{Proof of \fullref{lem:projection}}
\end{figure}

Note that, by the minimality assumption on $y$, we must have $x_{n-1} \in -\phi$. We have seen that both
$\partial \psi$ and $\partial r_\psi(\phi)$ separate $x$ from $y$. Since $\partial r_\psi(\phi) \neq \partial
\phi$ and since $\partial \phi$ is the only wall which separates $x_{n-1}$ from $y$, it follows that the
restricted path $\gamma'\co  x_0, x_1, \dots, x_{n-1}$ crosses both $\partial \psi$ and $\partial r_\psi(\phi)$,
but it does not cross $\partial \phi$. On the other hand we have $(-\phi) \cap \psi \subset -r_\psi(\phi)$ by
(1), which implies that $\gamma'$ crosses first $\partial \psi$ and then $\partial r_\psi(\phi)$. In other
words, there exist $i < j \in \{1, \dots, n-1\}$ such that $\psi$ contains $x_{i-1}$ but not $x_i$ and
$r_\psi(\phi)$ contains $x_j$ but not $x_{j-1}$. Now, the path $x_0, \dots, x_{i-1} = r_\psi(x_i),
r_\psi(x_{i+1}), \dots,  r_\psi(x_{j-1}), r_\psi(x_j)$ is of length $j < n$ and joins $x = x_0$ to $r_\psi(x_j)
\in \phi$. Again, this contradicts the minimality assumption on $y$. This proves (i).

Note that $x_{n-1} \in (-\phi) \cap (-\psi) \subset -r_\phi(\psi)$ by (2). Since $\partial \phi$ is the only
wall which separates $y$ from $x_{n-1}$ and since $\partial \phi \neq \partial r_\phi(\psi)$, we deduce $y \not
\in r_\phi(\psi)$. Assume in order to obtain a contradiction that $x \in r_\phi(\psi)$. We may then apply (i) to
the half-space $r_\phi(\psi)$, which yields $\phi \cap r_\phi(\psi) \subset \psi$. Transforming by $r_\phi$
yields $(-\phi) \cap \psi \subset r_\phi(\psi)$. By (i), we have also $\phi \cap \psi \subset r_\phi(\psi)$ and
we deduce $\psi = (\phi \cap \psi) \cup ((-\phi) \cap \psi) \subset r_\phi(\psi)$, whence $\psi = r_\phi(\psi)$.
This means that $r_\phi$ and $r_\psi$ commute, and we have seen above that it is not the case. This shows that
$x,y \not \in r_\phi(\psi)$, whence (ii).

It remains to prove (iii). Let thus $\alpha$ be a half-space containing $x$ and entirely contained in
$\psi$, and suppose that $\partial \alpha$ meets $\partial \phi$. We must prove that $\partial \alpha$ meets
$\partial r_\phi(\psi)$. Let $m$ be the midpoint of the edge joining $y$ to $x_{n-1}$; thus $m$ lies on
$\partial \phi$. By hypothesis, the half-space $\alpha$ contains $x$ but not $y$; therefore, the geodesic path
joining $x$ to $m$ must meet the wall $\partial \alpha$ in a point $a$. Let also $a'$ be a point lying on
$\partial \phi \cap \partial \alpha$. Consider the geodesic triangle whose vertices are $m, a, a'$. Since
$\partial \psi$ separates $m$ from $a'$ and since $\partial \psi \cap \partial \phi = \partial \phi \cap
\partial r_\phi(\psi)$, it follows that $\partial r_\phi(\psi)$ separates $m$ from $a'$. By (ii), it follows
that $\partial r_\phi(\psi)$ does not separate $m$ from $a$. Therefore, the wall $\partial r_\phi(\psi)$
separates $a$ from $a'$. Since both $a$ and $a'$ lie on $\partial \alpha$ and since walls are convex, we deduce
that $\partial \alpha$ meets $\partial r_\phi(\psi)$. This finishes the proof.
\end{proof}

\subsection{Existence of pairs of parallel walls}

\begin{lem}\label{lem:NR}
Suppose that $S$ is finite. Then there exists a constant $N=N(W,S)$ such that any set of more than $N$ walls
contains a pair of parallel walls.
\end{lem}
\begin{proof}
See \cite[Lemma 3]{NR03}.
\end{proof}

\subsection[Proof of \ref{thm:ParallelWalls}]{Proof of \fullref{thm:ParallelWalls}}

Note that $\mathscr{A}$ is quasi-isometric to its $1$--skeleton, which is the Cayley graph of $(W,S)$. Therefore,
if suffices to show the existence of a constant $Q$ such that any pair of parallel walls, at
\emph{combinatorial} distance greater than $Q$ from one another, is separated by another wall. This is what we
prove now; for convenience we denote by $V(\mathscr{A})$ the vertex set of $\mathscr{A}$ and by $d$ the
combinatorial distance on $V(\mathscr{A})$.

By \fullref{lem:NR} combined with Ramsey's theorem, there exists a constant $Q$ such that any set of more than
$Q$ walls contains a subset of $4l+1$ pairwise parallel walls, where $l= \inf\{ n \in \N \; | \; n
>L, n \geq 8\}$ and $L$ is the constant of \fullref{thm:ladder}.

Let $\alpha, \beta \in \Phi(\mathscr{A})$ be half-spaces such that $\alpha \cap \beta = \varnothing$ and let
$$d(\alpha, \beta) = \min \{d(x, y) \; | \; x, y \in V(\mathscr{A}), x \in \alpha, y \in \beta \}.$$
Assume that $d(\alpha, \beta) > Q$ and let $x, y \in V(\mathscr{A})$, $x \in \alpha$, $y \in \beta$ such that
$d(x, y) = d(\alpha, \beta)$. We must prove that there is a wall $m$ which separates $\partial \alpha$ from
$\partial \beta$; equivalently $m$ must separate $x$ from $y$, and be parallel to, but distinct from, $\partial
\alpha$ and $\partial \beta$.

The set $\mathscr{M}(x, y)$ consisting of all walls which separate $x$ from $y$ is of cardinality $d(x, y)$.
Thus it contains a subset $M$ of cardinality $4l+1$ consisting of pairwise parallel walls. If $\partial \alpha$
and $\partial \beta$ both belong to $M$ then we are done. Hence we may assume without loss of generality that
$\partial \beta \not \in M$. There are two cases to consider: either $\partial \alpha \not \in M$ or $\partial
\alpha \in M$. If $\partial \alpha$ is parallel to all elements of $M$, we may replace $M$ by $M \cup \{\partial
\alpha\}$ and we are reduced to the case when $\partial \alpha \in M$.

\begin{figure}[ht!]
\label{fig:ParallelWalls}
\begin{center}
\includegraphics[scale=0.48]{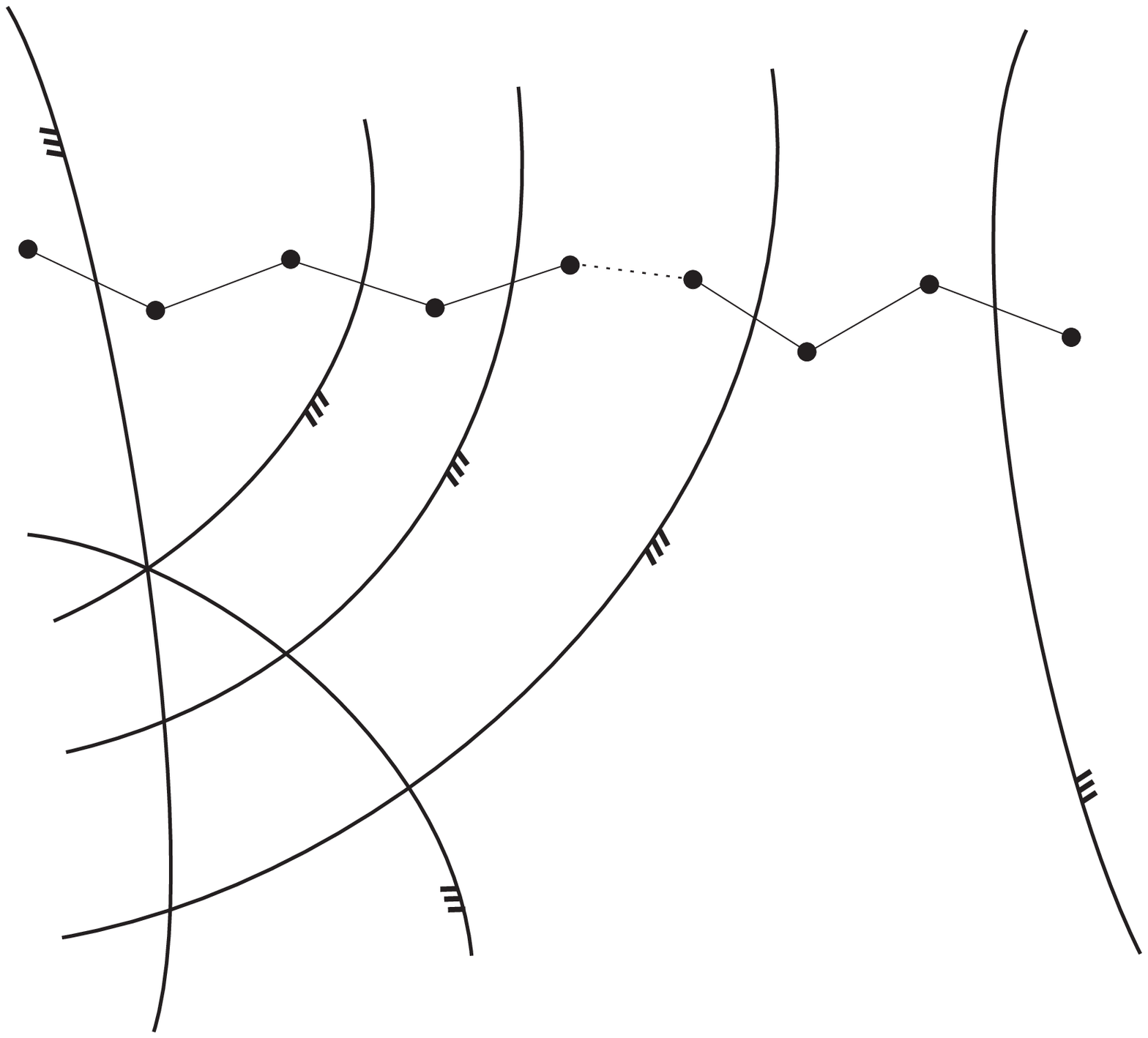}
\unitlength=.8pt\small
\begin{picture}(0,0)
\put(-218,25){$r_\alpha(\phi_{2l})$}%
\put(-290,185){$x$}%
\put(-26,187){$y$}%
\put(-96,155){$y'$}%
\put(-293,220){$\alpha$}%
\put(-13,60){$\beta$}%
\put(-123,115){$\phi_0$}%
\put(-173,135){$\phi_{2l-1}$}%
\put(-208,150){$\phi_{2l}$}%
\end{picture}
\end{center}
\caption{Proof of \fullref{thm:ParallelWalls}}
\end{figure}

We henceforth assume that $\partial \alpha \not \in M$ and that some element of $M$ meets $\partial
\alpha$. In fact, we may assume that each element of $M$ meets $\partial \alpha$ or $\partial \beta$, otherwise
we are done. Up to exchanging $\alpha$ and $\beta$, we may therefore assume that at least $2l+1$ elements of $M$
meet $\partial \alpha$. Let $\phi_0 \subsetneq \phi_1 \subsetneq \dots \subsetneq \phi_{2l}$ be half-spaces
containing $y$ but not $x$, such that $\partial \phi_i \in M$ meets $\partial \alpha$ for each $i$, see
\fullref{fig:ParallelWalls}. We have $\varnothing \neq \partial \alpha \cap \partial \phi_{2l} = \partial
\alpha \cap
\partial r_\alpha(\phi_{2l}) \subset
\partial \phi_{2l}$. We apply \fullref{lem:projection}(iii) with the half-spaces $\alpha$ and $\phi_{2l}$: this
shows that $\partial r_\alpha(\phi_{2l})$ meets $\partial \phi_i$ for each $i = 0, 1, \dots, 2l$. Therefore, the
hypotheses of \fullref{thm:ladder} are satisfied; this shows that $\la r_\alpha, r_{\phi_i} \; | \; i= 0, 1,
\dots, 2l \ra$ is a Euclidean triangle subgroup. Hence, we may apply \fullref{lem:nestedEucl}(ii) which yields
a set $M'$ of $l$ pairwise parallel walls which separate $\partial \alpha$ from any point on $\partial \phi_0
\cap
\partial r_\alpha(\phi_{2l})$; furthermore the group $W(M' \cup \{\partial \alpha\})$ is infinite dihedral.

We now show that $M' \subset \mathscr{M}(x, y)$. To this end, let $y'$ be the last vertex belonging to $\phi_0$,
on a geodesic path in $V(\mathscr{A})$ joining $y$ to $x$, and let $x'$ be any point on $\partial \alpha \cap
\partial \phi_0$, see \fullref{fig:ParallelWalls}.
By \fullref{lem:projection}(ii), the points $x$ and $y'$ lie on the same
side of $\partial r_\alpha(\phi_{2l})$. Moreover, since $\partial \phi_{2l}$ separates $x$ from $x'$ and since
$\partial \alpha \cap \partial \phi_{2l} = \partial \alpha \cap
\partial r_\alpha( \phi_{2l})$, it follows that $\partial r_\alpha( \phi_{2l})$ separates $x$ from $x'$.
Therefore $\partial r_\alpha( \phi_{2l})$ also separates $x'$ from $y'$. Hence, the geodesic path joining $x'$
to the $\CAT0$ projection of $y'$ on $\partial \phi_0$ meets $\partial r_\alpha( \phi_{2l})$ in a point $z$. We
have seen above that the elements of $M'$ separate $z$ from $\partial \alpha$; in particular they separate $z$
from $x$ and, hence, they separate $y'$ from $x$ because any geodesic segment crosses every wall at most once.
Since $\mathscr{M}(x, y') \subset \mathscr{M}(x, y)$, it follows that $M' \subset \mathscr{M}(x, y)$ as desired.

Thus all elements of $M'$ are parallel to $\partial \alpha$ and separate $x$ from $y$. If any element of $M'$ is
parallel to $\partial \beta$, then we are done. We henceforth assume that all elements of $M'$ meet $\partial
\beta$. Since $M'$ is of cardinality $l \geq 8$ and since $W(M' \cup \{\partial \alpha\})$ is infinite dihedral,
\fullref{lem:Rk3} implies that $\partial \alpha$ meets $\partial \beta$, which is absurd.

It remains to deal with the case when $\partial \alpha \in M$. If any element of $M$ distinct from
$\partial \alpha$ is parallel to $\beta$, then we are done. Thus we may assume that all other elements of $M$
meet $\partial \beta$. Repeating the same arguments as above appealing to \fullref{thm:ladder} and
\fullref{lem:nestedEucl}(ii), we obtain a set $M'' \subset \mathscr{M}(x, y)$, of cardinality $2l$, consisting
of pairwise parallel walls, which all separate $\partial \beta$ from $x$ and such that $W(M'' \cup \{\partial
\beta\})$ is infinite dihedral. Once again, if any element of $M''$ is parallel to $\partial \alpha$ then we are
done. Otherwise, we deduce from \fullref{lem:Rk3} that $\partial \beta$ meets $\partial \alpha$, a
contradiction. \qed

\section{The parabolic closure of a dihedral reflection subgroup}\label{sect:chainsII}

A basic fact on Coxeter groups is that any intersection of parabolic subgroups is itself a parabolic subgroup.
This allows to define the \textit{parabolic closure} $\Pc(R)$ of a subset $R \subset W$: it is the smallest
parabolic subgroup of $W$ containing $R$.

Given two roots $\alpha, \beta \in \Phi(\mathscr{E})$ such that $\alpha \neq \pm \beta$, it a well-known
consequence of \fullref{lem:LinkRootBases} that the dihedral group $\la r_\alpha, r_\beta \ra$ is infinite if
and only if $|(\alpha, \beta)| \geq 1$. The main result of this section is the following:

\begin{prop}\label{prop:AffinePairs}
Assume that $\mathscr{E}$ is the standard root basis. Let $\alpha, \beta \in \Phi(\mathscr{E})$ be such that
$\alpha \neq \pm \beta$. Then there exists a constant $\varepsilon >0$, depending only on the Coxeter system
$(W,S)$, such that the following assertions hold:
\begin{itemize}
\item[(i)] If $|(\alpha, \beta)| = 1$, then the parabolic closure of $\la r_\alpha, r_\beta \ra$ is of
irreducible affine type; in particular, it is virtually abelian.

\item[(ii)] If $|(\alpha, \beta)| > 1$, then $|(\alpha, \beta)| \geq 1 + \varepsilon$ and the parabolic closure
of $\la r_\alpha, r_\beta \ra$ is not of affine type; in particular, it contains a free subgroup of rank~$2$.
\end{itemize}
\end{prop}

The proof of \fullref{thm:chains} of the introduction, which is given in \fullref{sect:proof:chains}
below, is a combination of Propositions~\ref{prop:FiniteChain} and~\ref{prop:AffinePairs}.

\subsection{Orbits of pairs of walls}

Since $W$ has finitely many conjugacy classes of reflections, it follows that $W$ has finitely many orbits of
walls in $\mathscr{A}$. More generally, we have the following:

\begin{lem}\label{lem:OrbitsPairs}
For any $r \in \R_+$, the group $W$ has finitely many orbits of pairs of walls $\{m, m'\}$ such that $d(m, m')
\leq r$.
\end{lem}
\begin{proof}
If two walls meet, then the corresponding reflections generate a finite subgroup. Since $W$ has finitely many
conjugacy classes of finite subgroups, it follows immediately that $W$ has finitely many orbits of pairs of
intersecting walls. Therefore, it suffices to consider orbits of pairs of parallel walls, or equivalently, pairs
of disjoint half-spaces.

As in the proof of \fullref{thm:ladder}, we will consider only the combinatorial distance on
$V(\mathscr{A})$, which we denote again by $d$. To every pair of disjoint half-spaces $\{\alpha, \beta\}$, we
associate an oriented combinatorial path $\gamma(\alpha; \beta)$ in the $1$--skeleton $\mathscr{A}_0$ as follows:
Choose $x, y \in V(\mathscr{A})$ such that $x \in \alpha$, $y \in \beta$ and $d(x, y)$. Since $\alpha$ and
$\beta$ are disjoint, the walls $\partial \alpha$ and $\partial \beta$ both belong to $\mathscr{M}(x, y)$ and,
hence, there exists a geodesic path joining $x$ to $y$ in $\mathscr{A}_0$, whose first edge crosses $\partial
\alpha$ and whose last edge crosses $\partial \beta$. We define $\gamma(\alpha; \beta)$ to be any such geodesic
path. Clearly the length of $\gamma(\alpha; \beta)$ equals $d(\alpha, \beta)$ and moreover, if $\{\alpha',
\beta'\}$ is a pair of disjoint roots such that $w.\gamma(\alpha'; \beta') = \gamma(\alpha; \beta)$ for some $w
\in W$, then $w.\alpha' =  \alpha$ and $w.\beta' = \beta$. Therefore, in order to prove the desired finiteness
property, it suffices to show that $W$ has finitely many orbits on oriented combinatorial paths of length~$\leq
r$. This is true because $W$ is transitive on $V(\mathscr{A})$ and $\mathscr{A}$ is locally finite.
\end{proof}

\subsection{On the parabolic closure}

By definition of the Davis complex $\mathscr{A}$, the stabilizer of any point of $\mathscr{A}$ in $W$ is a
finite parabolic subgroup. Since any finite subgroup of $W$ fixes a point of $\mathscr{A}$, it follows that any
finite subgroup is contained in a finite parabolic subgroup; in other words, if a subgroup of $W$ is finite,
then so is its parabolic closure. This is well-known. A fundamental result, due to D~Krammer, is that this is
also true for affine reflection subgroups:

\begin{prop}\label{prop:PcAffine}
Given any reflection subgroup $R$, if $R$ is of irreducible affine type and rank~$\geq 3$, then so is its
parabolic closure.
\end{prop}
\begin{proof}
Notice first that if a reflection subgroup is of irreducible type, then so is its parabolic closure. Therefore,
the desired assertion follows from \cite[Theorem~3.3]{CM05}, which builds upon the results of
\cite[Chapter~6]{Kr94}.
\end{proof}

The following lemma provides useful criterions ensuring that a given reflection belongs to a certain parabolic
closure:

\begin{lem}\label{lem:Pc}
Let $M$ be a set of cardinality at least~$2$, consisting of pairwise parallel walls and let $m$ be any wall.
Assume that any of the following conditions hold:
\begin{itemize}
\item[(1)] $m$ separates two elements of $M$.

\item[(2)] $W(M \cup \{m\})$ is infinite dihedral.

\item[(3)] $W(M \cup\{m\})$ is a Euclidean triangle subgroup.
\end{itemize}
Then the reflection $r_m$  belongs to the parabolic closure $\Pc(W(M))$.
\end{lem}
\begin{proof}
Point (1) follows from the same convexity arguments as in the proof of \fullref{lem:nestedEucl}(i).

(2)\qua The hypotheses imply that $W(M)$ is infinite dihedral. Let $t$ be a nontrivial translation of $W(M)$ and
let $\mu \in M$. Since $W(M \cup \{m\})$ is infinite dihedral, there exists a nonzero integer $n \in \Z$ such
that $m$ separates $\mu$ from $t^n.\mu$. By (1), this yields $r_m \in \Pc(W(\mu, t^n.\mu))$. On the other hand,
we have clearly $W(\mu, t^n.\mu) \subset W(M)$, whence $\Pc(W(\mu, t^n.\mu)) \subset \Pc(W(M))$.

(3)\qua By \fullref{prop:PcAffine}, the parabolic subgroup $\Pc(W(M \cup\{m\}))$ is of irreducible affine
type. On the other hand $W(M) \subset W(M \cup\{m\})$ and hence $\Pc(W(M)) \subset \Pc(W(M \cup\{m\}))$. Since
any proper parabolic subgroup of a parabolic of irreducible affine type is finite while $W(M)$ is infinite, we
deduce that $\Pc(W(M)) = \Pc(W(M \cup\{m\}))$.
\end{proof}

\subsection[Proof of \ref{prop:AffinePairs}]{Proof of \fullref{prop:AffinePairs}}

In order to simplify notation, we assume that a base chamber $C \subset \mathscr{A}$ has been fixed and we
identify $\Phi(\mathscr{E})$ with $\Phi(\mathscr{A})$ be means of the bijection $\zeta_C$ of
\fullref{lem:LinkRootBases}(i), and we note $\Phi = \Phi(\mathscr{E}) = \Phi(\mathscr{A})$.

\begin{proof}[Proof of \fullref{prop:AffinePairs}(i)]
Let $\alpha, \beta \in \Phi$ be roots such that $\alpha \neq \pm \beta$ and $|(\alpha, \beta)| = 1$. We must
prove that $\Pc(r_\alpha, r_\beta)$ is of irreducible affine type. Since $\la r_\alpha, r_\beta \ra$ is infinite
and since every proper parabolic subgroup of a Coxeter group of irreducible affine type is finite, it suffices
to prove that $r_\alpha$ and $r_\beta$ are contained in a common parabolic subgroup of irreducible affine type.

Up to replacing $\alpha$ or $\beta$ by its opposite, we may assume without loss of generality that $\alpha
\subset \beta$, whence $(\alpha, \beta) =1$ by \fullref{lem:LinkRootBases}(iii).

Let $x, y \in V(\mathscr{A})$ be vertices such that $x \in \alpha$, $y \in -\beta$ and $d(x, y) = d(\alpha,
-\beta)$. Let $\Phi(x, y)$ be the set of all half-spaces containing $x$ but not $y$.  By
\fullref{lem:Phi(alpha;beta)}, the set $\Phi(\alpha; \beta)$ is a chain. We denote it by
$$\alpha = \alpha_0 \subsetneq \alpha_1 \subsetneq \dots \subsetneq \alpha_k = \beta.$$
Moreover, \fullref{lem:Phi(alpha;beta)} implies that the group $W(\Phi(\alpha; \beta)) = \la r_{\alpha_i} \; |
\; i=0, \dots, k \ra$ is infinite dihedral, generated by $\{r_{\alpha_0}, r_{\alpha_1}\}$.

Let $t = r_\alpha r_\beta$. A straightforward computation shows that $(\alpha, t^n. \beta) = 1$ for all $n \in
\Z$. Moreover, we have $r_\beta \in \Pc(\la r_\alpha, r_{t^n.\beta}\ra)$ for all $n \neq 0$ by
\fullref{lem:Pc}. Therefore, up to replacing $\beta$ by $t^n.\beta$ with $n$ sufficiently large, we may assume
without loss of generality that $k \geq 8$.

Assume first that $\Phi(\alpha; \beta) = \Phi(x, y)$. Then, considering a geodesic path joining $x$ to $y$ in
the $1$--skeleton of $\mathscr{A}$, we see that there is a combinatorial path of length~$2$ crossing successively
the walls $\partial \alpha_0$ and $\partial \alpha_1$. This means that the infinite dihedral group
$W(\Phi(\alpha; \beta)) =\la r_{\alpha_0}, r_{\alpha_1} \ra$ is a parabolic subgroup of rank~$2$. Thus it is of
irreducible affine type and it contains $r_\beta$ since $\beta \in \Phi(\alpha; \beta)$. Hence we are done in
this case.

Assume now that $\Phi(\alpha; \beta)$ is properly contained in $\Phi(x, y)$ and consider an element $\gamma \in
\Phi(x, y) \backslash \Phi(\alpha; \beta)$. By the definition of $\Phi(\alpha; \beta)$, it follows that
$\partial \gamma$ must meet $\partial \alpha$ or $\partial \beta$. Without loss of generality, we assume that
$\partial \gamma$ meets $\partial \beta$. Let
$$\Phi_0(\gamma) = \{\phi \in \Phi(\alpha; \beta) \; | \; \partial \phi \text{ meets } \partial
\gamma \}.$$ Since walls are convex, it follows that if $\alpha_i, \alpha_{i'} \in \Phi_0(\gamma)$, then
$\alpha_j \in \Phi_0(\gamma)$ for all $i \leq j \leq i'$.

Note that for each $\phi \in \Phi(x, y)$, the reflections $r_\beta$ and $r_\phi$ do not commute, otherwise we
would have $x \in \phi = r_\beta(\phi)$, in contradiction with \fullref{lem:projection}(ii).

Suppose that $|\Phi_0(\gamma)| > 7$; this happens whenever $\alpha \in \Phi_0(\gamma)$. Since $W(\Phi(\alpha;
\beta))$ is infinite dihedral, this implies by \fullref{lem:Rk3} that either $r_\gamma$ centralizes the group
$W(\Phi(\alpha; \beta))$ or that $W(\Phi(\alpha; \beta))$ is contained in a Euclidean triangle subgroup. We have
just seen $r_\gamma$ does not centralize $r_\beta \in W(\Phi(\alpha; \beta))$, so $W(\Phi(\alpha; \beta))$ is
contained in a Euclidean triangle subgroup. By \fullref{prop:PcAffine}, a Euclidean triangle subgroup is
contained in a parabolic subgroup of irreducible affine type, hence we are done in this case.

Similarly, if the wall $\partial r_\beta(\gamma)$ meets the boundary of more than $7$~elements of $\Phi(\alpha;
\beta)$, it follows also that $W(\Phi(\alpha; \beta))$ is contained in a Euclidean triangle group, because
$r_{r_\beta(\gamma)}= r_\beta r_\gamma r_\beta$ does not commute with $r_\beta$ since $r_\gamma$ does not. Hence
we are done in this case as well.

It remains to consider the case when $\partial \alpha$ meets neither $\partial \gamma$ nor $\partial
r_\beta(\gamma)$. This yields $\alpha \subset \gamma$ and $\alpha \subset -r_\beta(\gamma)$, whence $(\alpha,
\gamma) \geq 1$ and $(\alpha, r_\beta(\gamma)) \leq -1$ by \fullref{lem:LinkRootBases}(iii). Therefore we
obtain
$$2 \leq 1+(\alpha, \gamma) \leq 2(\alpha, \beta)(\beta, \gamma).$$ By \fullref{lem:LinkRootBases}(iv)
and \fullref{lem:projection}(i), we have $(\beta, \gamma) > 0$. Thus we obtain
$$(\alpha, \beta) \geq \frac{1}{(\beta, \gamma)}.$$
On the other hand $( \beta,  \gamma) < 1$ by \fullref{lem:LinkRootBases}(ii), because $\partial \beta$ meets
$\partial \gamma$. This contradicts the fact that $(\alpha, \beta)=1$, thereby proving that this last case does
not occur.
\end{proof}

\begin{proof}[Proof of \fullref{prop:AffinePairs}(ii)]
If $(W,S)$ is of irreducible affine type, then Tits bilinear form is positive semi-definite \cite[Chapter~VI, Section~4.3,
Theorem~4]{Bo81}. Therefore, for all $\alpha, \beta \in \Phi$ we have $|(\alpha, \beta)| \leq \sqrt{(\alpha, \alpha)
(\beta, \beta)} = 1$ by Cauchy--Schwarz. This shows that if $(W,S)$ is arbitrary and if $\Pc(r_\alpha, r_\beta)$
is of irreducible affine type, then $|(\alpha, \beta)| \leq 1$.

It remains to prove that $\varepsilon > 0$, where
$$\varepsilon = -1 +\inf\{ (\alpha, \beta)  \; | \; \alpha, \beta \in \Phi, \ (\alpha, \beta) > 1 \}. $$
To this end, we define $P_0$ to be the set of all ordered pairs $(\alpha ; \beta) \in \Phi \times \Phi$ such
that
\begin{itemize}
\item[(1)] $\alpha \subset \beta$,

\item[(2)] $(\alpha, \beta) > 1$,

\item[(3)] For all $\gamma \in \Phi$ such that $\alpha \subset \gamma \subset \beta$, one has $\gamma \in
\{\alpha, \beta\}$.
\end{itemize}
We also define $P_1$ to be the set of all ordered pairs $(\alpha ; \beta) \in \Phi \times \Phi$ satisfying (1)
and (2) but not (3). In view of \fullref{lem:LinkRootBases}(iii), we have
$$1+ \varepsilon = \inf\{ (\alpha, \beta)  \; | \; (\alpha ; \beta) \in P_0 \cup P_1 \}. $$
In other words, we have $\varepsilon = \min \{\varepsilon_0, \varepsilon_1\}$, where $1 + \varepsilon_i = \inf\{
(\alpha, \beta)  \; | \; (\alpha ; \beta) \in P_i \}$ for $i = 0,1$.

Condition (3) means that the walls $\partial \alpha$ and $\partial \beta$, which are parallel, are not separated
by any wall. By \fullref{thm:ParallelWalls}, this implies that the distance from $\partial \alpha$ to
$\partial \beta$ is at most $Q$. Therefore, the group $W$ has finitely many orbits in $P_0$ by
\fullref{lem:OrbitsPairs}. In particular the set $\{ (\alpha, \beta)  \; | \; (\alpha ; \beta) \in P_0 \}$ is
finite and, hence, we have $\varepsilon_0 > 0$.

Let now $(\alpha ; \beta) \in P_1$. The set $\Phi(\alpha;\beta) = \{ \phi \in \Phi \; | \; \alpha \subset \phi
\subset \beta\}$ is finite. Since $(\alpha ; \beta) \not \in P_0$, there exists $\gamma \in \Phi(\alpha;\beta)$
distinct from $\alpha$ and $\beta$. Among all such $\gamma$'s, we choose one which  is minimal for $\subset$. In
particular, the pair $(\alpha ; \gamma)$ satisfies condition (3). Since $\alpha \subset \gamma$ we have
$(\alpha, \gamma) \geq 1$ by \fullref{lem:LinkRootBases}(iii). There are two cases.

First, assume that $(\alpha, \gamma) > 1$. In that case we have $(\alpha ; \gamma) \in P_0$ whence $(\alpha ,
\gamma) \geq 1+ \varepsilon_0$. By \fullref{lem:OrbitsPairs}(i), this yields $(\alpha, \beta) \geq 1+
\varepsilon_0$.

Assume now that $(\alpha, \gamma) = 1$. Choose $\gamma' \in \Phi(\alpha;\beta)$ such that $(\alpha, \gamma') =
1$ and such that $\gamma'$ is maximal (for $\subset$) with respect to these properties. By
\fullref{lem:Phi(alpha;beta)}, the set $\Phi(\alpha ; \gamma')$, which is contained in $\Phi(\alpha ; \beta)$,
is a chain. Let $\gamma''$ be the maximal element of $\Phi(\alpha ; \gamma') \backslash \{\gamma'\}$ and
consider $r_{\gamma'}(\gamma'')$. Using \fullref{lem:Phi(alpha;beta)}, it is immediate to compute that
$(\alpha, r_{\gamma'}(\gamma'')) = -1$; in particular, by maximality of $\gamma'$ we have
$-r_{\gamma'}(\gamma'') \not \in \Phi(\alpha ; \beta)$, or in other words $-r_{\gamma'}(\gamma'') \not \subset
\beta$. Similarly, if $\beta \subset -r_{\gamma'}(\gamma'')$, then \fullref{lem:OrderedTriples}(i) yields
$(\alpha, \beta) \leq (\alpha, -r_{\gamma'}(\gamma'')) =1$, contradicting (2). Therefore,
\fullref{lem:OrderedTriples}(iii) implies that $(r_{\gamma'}(\gamma''), \beta) > -1$. Using
\fullref{lem:OrderedTriples}(i) and~(ii), this yields successively
$$(\alpha, \beta) \  \geq \ (\gamma'', \beta) \  \geq \ 2(\gamma', \beta) - \kappa \ \geq \ 2 - \kappa.$$

The preceding two paragraphs show that for all $(\alpha ; \beta) \in P_1$, we have $(\alpha, \beta) \geq \max
\{1+\varepsilon_0, 2- \kappa\}$, whence $\varepsilon_1 = \min \{ \varepsilon_0, 1 - \kappa\}$. Recall from
\fullref{sect:kappa} that $\kappa <1$; furthermore, we have seen above that $\varepsilon_0 > 0$. This shows
finally that $\varepsilon >0$, as desired.
\end{proof}

\subsection[Proof of \ref{thm:chains}]{Proof of \fullref{thm:chains}}\label{sect:proof:chains}

Let $\alpha_0 \subsetneq \alpha_1 \subsetneq \dots \subsetneq \alpha_n$ be a sequence of roots. If $(\alpha_0,
\alpha_n)=1$, then one is in the alternative~(2) of the theorem by Lemmas~\ref{lem:Phi(alpha;beta)},
\ref{lem:Pc} and \fullref{prop:AffinePairs}(i). Thus we may assume that $(\alpha_0, \alpha_n)>1$ and we
must show that $(\alpha_0, \alpha_n)$ is bounded from below by a non-decreasing function of $n$.

In order to estimate $(\alpha_0, \alpha_n)$, we choose a maximally convex chain $\beta_0 \subsetneq \dots
\subsetneq \beta_m$ such that $\beta_0 = \alpha_0$ and $\beta_m = \alpha_n$. By definition, we have $m \geq n$.
Let $j = \min \{ i \geq 0 \; | \; (\beta_0, \beta_i) > 1\}$. Thus $j > 0$. By
\fullref{prop:FiniteChain}, we have $(\beta_0, \beta_m) > \max \{j ( 1 - \kappa), \  1 + \frac{m}{2
j}\varepsilon\}$, where $\varepsilon$ is the constant of \fullref{prop:AffinePairs}(ii). One easily
computes that
$$ j ( 1 - \kappa) \geq   1 + \frac{m}{2
j}\varepsilon \hspace{1cm}\Leftrightarrow \hspace{1cm} j \geq \frac{-1+\sqrt{1 + 2 \varepsilon(1-\kappa) m}}{1 -
\kappa}.$$ Therefore, one has
$$(\beta_0, \beta_m) > \sqrt{1 + 2 \varepsilon(1-\kappa) m}-1
\hspace{1cm}\text{if} \hspace{1cm} j \geq \frac{-1+\sqrt{1 + 2 \varepsilon(1-\kappa) m}}{1 - \kappa}$$ and
$$ (\beta_0, \beta_m) > 1 + \frac{\varepsilon (1-\kappa )m}{-2+2\sqrt{1 + 2 \varepsilon(1-\kappa) m}}
\hspace{1cm}\text{otherwise.}$$ Set $r_0 = 1$ and for all positive integer $k$, define $$r_k = \min \Big\{ -1
+\sqrt{1 + 2 \varepsilon(1-\kappa) k} , \ 1 + \frac{\varepsilon (1-\kappa )k}{-2 + 2\sqrt{1 + 2
\varepsilon(1-\kappa) k}} \Big\}.$$ Note that $(r_k)_{k \geq 1}$ is a non-decreasing sequence with $r_1 > 1$,
which tends to $+\infty$ with $k$. We have seen above that $(\beta_0, \beta_m) \geq r_m$. Since $\beta_0 =
\alpha_0$ and $\beta_m = \alpha_n$ and $m \geq n$, this yields $(\alpha_0, \alpha_n) \geq r_n$, and we are in
the alternative~(1) of the theorem. This finishes the proof. \qed

\section{Conjugacy of $2$--spherical subgroups} \label{sect:FiniteConjug}
\vspace{-5pt}

The purpose of this section is to prove \fullref{thm:FiniteConjug:2sph} and its corollary. As mentioned in
the introduction, this first requires to study in details  some aspects of the $\CAT 0$ cube complex
$\mathscr{X}$ of G~Niblo and L~Reeves; see \fullref{prop:CubicalChamber} and \fullref{thm:cubes}
below for specific statements.
\vspace{-5pt}

\subsection{Pairwise intersecting walls}\label{sect:M(c)}
\vspace{-5pt}

Any $n$--dimensional cube $c$ of a $\CAT 0$ cube complex determines a unique $n$--tuple of walls, denoted by
$M(c)$, consisting of the walls which contain the center of that cube. The following basic fact is an important
property:
\vspace{-5pt}

\begin{lem}\label{lem:Pairwise}
Let $c, c'$ be cubes of a $\CAT 0$ cube complex such that every wall in $M(c)$ meets every wall in $M(c')$. Let
$d$ be the combinatorial distance from $c$ to $c'$. Then there exists a cube $c''$ at combinatorial distance at
most $d$ from $c$, such that $M(c'') = M(c) \cup M(c')$.
\end{lem}
\vspace{-5pt}
\begin{proof}
See \cite[Theorem 4.14]{Sageev}. The estimate of the distances follows from an easy induction.
\end{proof}
\vspace{-5pt}

The only cube complex we will need to consider here is the $\CAT 0$ complex $\mathscr{X}$. We have seen in
\fullref{sect:NR} that the walls and hyperplanes of $\mathscr{X}$ are in canonical one-to-one correspondence
with those of $\mathscr{A}$. In order to simplify notation, we identify $\mathscr{M(A)}$ with $\mathscr{M(X)}$
(resp.\ $\Phi(\mathscr{A})$ with $\Phi(\mathscr{X})$). A fundamental difference between $\mathscr{A}$ and
$\mathscr{X}$ is that a reflection of $W$ fixes a wall of $\mathscr{A}$ pointwise, while in $\mathscr{X}$ a
reflection acts mostly non-trivially on the wall it stabilizes. This explains why the property of
\fullref{lem:Pairwise} cannot hold in $\mathscr{A}$, unless $(W,S)$ is right-angled in which case $\mathscr{A}
= \mathscr{X}$. However, pairs of walls behave always similarly in $\mathscr{A}$ and in $\mathscr{X}$:
\vspace{-5pt}

\begin{lem}\label{lem:PairsHalfspaces}
Two walls meet in $\mathscr{A}$ if and only if they meet in $\mathscr{X}$. Two half-spaces are nested (resp.\
have empty intersection) in $\mathscr{A}$ if only if they are nested (resp.\ have empty intersection) in
$\mathscr{X}$.
\end{lem}
\vspace{-5pt}
\begin{proof}
This is a straightforward consequence of the construction of $\mathscr{X}$, see \fullref{sect:NR} and
\cite{NR03} for more details.
\end{proof}

Note that the corresponding statement fails for triples of walls. Indeed, every triple of pairwise intersecting
walls has nonempty intersection in $\mathscr{X}$ by \fullref{lem:Pairwise}, but this is obviously false in
$\mathscr{A}$: any Euclidean or hyperbolic triangle group provides a counterexample.

As above, for any set of walls $M$, we denote by $W(M)$ the subgroup of $W$ generated by all reflections through
elements of $M$.

\begin{lem}\label{lem:W(M):2sph}
Let $M$ be a set of pairwise intersecting walls. Then $W(M)$ is of $2$--spherical type.
\end{lem}
\begin{proof}
By \fullref{lem:PairsHalfspaces}, if two walls meet in $\mathscr{X}$, then they meet in $\mathscr{A}$ and,
hence, the corresponding reflections generate a finite subgroup. This shows that $W(M)$ is $2$--spherical in the
sense of \fullref{sect:2sph}. We have seen in this this section that $2$--spherical Coxeter groups are
precisely those Coxeter groups of $2$--spherical type, in the usual sense. This means that $W(M)$ is of
$2$--spherical type as a reflection subgroup.
\end{proof}

\subsection{On standard parabolic subgroups}

Recall that a \textit{standard parabolic subgroup} of $W$ is a subgroup generated by some subset of $S$, and a
subgroup is \textit{parabolic} if it is conjugate to a standard parabolic subgroup. In some situations, it is
useful to keep track of an element of $W$ which conjugates a given parabolic subgroup to a standard parabolic
one. This motivates the following definition: given a vertex $v_0$ (or a chamber) of $\mathscr{A}$ and a set
$M_0$ of walls such that each element of $M_0$ separates $v_0$ from a neighboring vertex, we say that the
parabolic subgroup $W(M_0)$ is \textit{standard with respect to} $v_0$. Thus a standard parabolic subgroup is a
parabolic subgroup which is standard with respect to the base chamber of $\mathscr{A}$ (ie, the chamber which
corresponds to the identity $1 \in W$ in the Cayley graph).

\begin{lem}\label{lem:Pc:standard}
Let $M$ be a set of pairwise parallel walls, of cardinality~$>7$, such that the parabolic closure $P =
\Pc(W(M))$ is of affine type. Let $v_0$ be a vertex of $\mathscr{A}$. Assume that, given any wall $m$, if $m$
separates $v_0$ from some wall in $M$, then $m \in M$. Then the parabolic subgroup $P$ is standard with respect
to $v_0$.
\end{lem}
\begin{proof}
Given a parabolic subgroup $W_0$ of $W$ and a chamber $c$ of $\mathscr{A}$ such that $W_0$ is standard with
respect to $c$, then the union of the $W_0$--orbit of $c$ is a closed convex subset of $\mathscr{A}$, which we
call a $W_0$--\emph{residue}. Clearly, the parabolic subgroup $W_0$ is standard with respect to a given chamber
if and only if this chamber is contained in some $W_0$--residue.

Let $c_0$ be the unique chamber of $\mathscr{A}$ containing $v_0$ and let $\rho_0$ be the $W_0$--residue at
minimal combinatorial distance from $c_0$. We must prove that $c_0 \subset \rho_0$. Assume the contrary in order
to obtain a contradiction. Let $c'$ be the chamber of $\rho_0$ at minimal combinatorial distance from $c_0$, let
$c''$ be a chamber adjacent to $c'$ and closer to $c_0$ than $c'$. Finally, let $m$ be the wall which separates
$c'$ from $c''$.

The reflection $r_m$ does not belong to $P$. Indeed, if $r_m \in P$, then $r_m$ would stabilize $\rho_0$ which
would imply that $c'' = r_m(c')$ be contained in $\rho_0$, in contradiction with the definition of $c'$. It
follows in particular that the wall $m$ separates $c_0$ from $\rho_0$: otherwise $m$ would separate two adjacent
chambers contained in $\rho_0$ and hence $r_m$ would swap these chambers, but, on the other hand, the only
element of $W$ swapping these chambers belongs to $P$.

Since $r_m \not \in P$ we have $m \not \in M$. Moreover, any wall $m' \in M$ meets $\rho_0$ because the
reflection $r_{m'} \in P$ stabilizes $\rho_0$ which is closed and convex. Therefore, the wall $m$ meets $m'$
otherwise $m$ would separate $v_0$ from $m'$ which is excluded by hypothesis. Thus $m$ meets each element of
$M$.

Since $M$ has at least $8$~elements, it follows from \fullref{lem:Rk3} that either $r_m$ centralizes $W(M)$ or
that $W(M \cup \{m\})$ is a Euclidean triangle subgroup. The second case is impossible by \fullref{lem:Pc}(3)
because $r_m \not \in P= \Pc(W(M))$. Thus $r_m$ centralizes $W(M)$ and, hence, normalizes the parabolic closure
$P$ of $W(M)$.

The set $r_m(\rho_0)$ is a $(r_m P r_m)$--residue and, hence, a $P$--residue by the preceding paragraph. Since
$r_m(\rho_0)$ contains the chamber $c'' = r_m(c')$, we obtain a contradiction with the minimality assumption we
made on $\rho_0$. This finishes the proof.\end{proof}

\subsection{The normalizer of an affine parabolic subgroup}

The following fact is well-known; it is more generally true for any infinite parabolic subgroup of irreducible
type.

\begin{lem}\label{lem:Normalizer}
Let $P \subset W$ be a parabolic subgroup of irreducible affine type. Then the normalizer of $P$ in $W$ splits
as a direct product: $N_W(P) = P \times C_W(P)$. In particular,  any reflection which normalizes $P$ either
belongs to $P$ or centralizes $P$.
\end{lem}
\begin{proof}
See \cite[Proposition~5.5]{Deodhar82}.
\end{proof}

\subsection{Free abelian normal subgroups in Coxeter groups}

The following statement of independent interest is a consequence of the work of Daan Krammer
\cite[Section 6.8]{Kr94} on free abelian subgroups of Coxeter groups:

\begin{lem}\label{lem:FreeAbNormal}
The group $W$ possesses a nontrivial free abelian normal subgroup if and only if the Coxeter diagram of $(W,S)$
has a connected component of irreducible affine type.
\end{lem}
\begin{proof}
If $(W,S)$ has a connected component of affine type, then the translation subgroup of the corresponding affine
parabolic subgroup is a nontrivial free abelian normal subgroup. Thus the `if' part is clear. Conversely, let
$H$ be a nontrivial free abelian normal subgroup of $W$. Let $W=W_1 \times \dots W_k$ be the decomposition of
$W$ into its direct components and let $\pr_i$ be the canonical projection of $W$ onto $W_i$. By Selberg's lemma
$W$ has a finite index torsion free subgroup. In particular $H$ has a finite index subgroup $H'$ such that
$\pr_i(H')$ is torsion free for all $i$. Since $H$ is free abelian, we may assume without loss of generality
that $H'$ is normalized by $W$. In particular $\pr_i(H')$ is a free abelian normal subgroup of $W_i$ for each
$i$. This shows that, in order to finish the proof, it suffices to show that an \emph{irreducible} Coxeter group
which possesses a nontrivial free abelian normal subgroup must be of affine type.

We assume henceforth that $(W,S)$ is irreducible, but not of affine type otherwise we are done. Since $H$ is
normal in $W$, so is its parabolic closure $\Pc(H)$. Since a parabolic subgroup is normal if and only if it is a
direct component, it follows that $\Pc(H) = W$ since $(W,S)$ is irreducible. This is true even after replacing
$H$ by a finite index subgroup normalized by $W$. Therefore, \cite[Theorem~6.8.3]{Kr94} implies that $H$ is of
rank one because $(W,S)$ is not of affine type. In particular the centralizer of $H$ in $W$ is of index at
most~$2$. On the other hand $H$ is of finite index in its centralizer by \cite[Corollary~6.3.10]{Kr94}. This
shows that $H$ is of finite index in $W$. Since $W$ is not of affine type, it follows that $W$ must be finite
and, hence, that $H$ is trivial. This is a contradiction.
\end{proof}

\subsection{The cubical chamber}

Recall from \fullref{sect:NR} that the Cayley graph of $(W,S)$ is equivariantly embedded in the $1$--skeleton
$\mathscr{X}^{(1)}$ of $\mathscr{X}$. We denote this subgraph by $\mathscr{X}_0$. Given a vertex $v \in
\mathscr{X}_0$, let $\Psi(v)$ be the set of those half-spaces which contain $v$ but not its neighbors in
$\mathscr{X}_0$. The set $\bigcap_{\psi \in \Psi(v)} \psi$, viewed as a subset of $\mathscr{X}$, is called the
\textit{cubical chamber} containing $v$. Given two points $x, y \in \mathscr{X}$ we denote by $\mathscr{M}(x,
y)$ the set of walls which separate $x$ from $y$.

\begin{prop}\label{prop:CubicalChamber}
Let $v_0$ be a vertex of $\mathscr{X}_0$ and let $x$ be a vertex of $\mathscr{X}$ belonging to the cubical
chamber containing $v_0$. Let $M \subset \mathscr{M}(v_0, x)$ be a set of pairwise parallel walls. There exists
a constant $K = K(W,S)$ such that if $M$ has more than $K$ elements, then $W(M)$ is infinite dihedral and its
parabolic closure $\Pc(W(M))$ is of irreducible affine type and rank~$\geq 3$.
\end{prop}
\begin{proof}
Define a constant $\lambda_{\max}$ as follows:
$$\lambda_{\max} = \sup \{ \lambda \in \R_+ \; | \; \lambda
  \text{ is a coefficient of some } \phi \in \Phi(\mathscr{E})^+_{\min}
  \text{ in the basis } \Pi\},$$
where
$$\Phi(\mathscr{E})^+_{\min} = \{ \alpha \in \Phi(\mathscr{E})^+ \; | \; \phi \in \Phi(\mathscr{E})^+, \phi
\subset \alpha \Rightarrow \phi = \alpha\}.$$
It follows from the parallel wall theorem \cite[Theorem~2.8]{BH94}
that the set $\Phi(\mathscr{E})^+_{\min}$ is finite. Therefore, the constant $\lambda_{\max}$ is a well-defined
positive real number.
\vspace{4pt}

Let also $r = |S|$ be the rank of $(W,S)$ and $\kappa$ be the constant defined in \fullref{sect:kappa}. We
will show that the desired constant $K$ can be defined as $K = \min \{ n \in \N \; | \; r_n > r \kappa
\lambda_{\max}\}$, where $(r_n)_{n \in \N}$ is the sequence of \fullref{thm:chains}.
\vspace{4pt}

Let $\Phi(x, v_0)$ be the set of those half-spaces which contain $x$ but not $v_0$. The set
$\Phi(x, v_0)$ contains a nested sequence $\phi_0 \subsetneq \phi_1 \subsetneq \dots \subsetneq
\phi_l$ of half-spaces such that $M \subset \{ \partial \phi_i  \; | \; i = 0, 1, \dots, l\}$. Without loss of
generality, we may -- and shall -- assume that $(\phi_i)_{i \leq l}$ is a maximal nested sequence contained in
$\Phi(x, v_0)$: since $\Phi(x, v_0)$ is finite, any nested sequence contained in $\Phi(x, v_0)$ can be completed
in order to obtain a maximal nested sequence. Note that if $(\phi_i)_{i \leq l}$ is maximal, then no wall
separates $v_0$ from $\partial \phi_l$, because if such a wall existed, we could lengthen the nested sequence
$(\phi_i)_{i \leq l}$ of one unit by adding the half-space containing $x$ and determined by this extra wall.
\vspace{4pt}

Let $\Psi$ be the set of those half-spaces which contain $v_0$ but not its neighbors in $\mathscr{X}_0$ and let
$\Psi_0 = \{ \psi \in \Psi \; | \; \phi_0 \not \subset \psi\}$. Note that every half-space contains a vertex of
$\mathscr{X}_0$. Furthermore, it follows from  the definition of $\Psi$ that the only vertex of $\mathscr{X}_0$
which is contained in $\bigcap_{\psi \in \Psi} \psi$ is $v_0$. Therefore, we deduce that for each $\phi \in
\Phi(x, v_0)$, there exists $\psi \in \Psi$ such that $\phi \not \subset \psi$. In particular, the set $\Psi_0$
is nonempty. Since the sequence $(\phi_i)_{i \leq l}$  is nested, we deduce that $\phi_i \not \subset \psi$ for
all $i=0, 1, \dots, l$ and all $\psi \in \Psi_0$.
\vspace{4pt}

Note that for all $i \in \{0,1, \dots,l \}$ and all $\psi \in \Psi_0$, the vertex $x$ is contained in $\phi_i
\cap \psi$ while the vertex $v_0$ is contained in $\psi$ but not in $\phi_i$. Since $\phi_i \not \subset \psi$
and since $v_0$ belongs to an edge crossed by $\partial \psi$, it then follows that the walls $\partial \phi_i$
and $\partial \psi$ meet for all $i \in \{0,1, \dots, l\}$ and all $\psi \in \Psi_0$.
\vspace{4pt}

 Let us choose as base chamber $C \subset \mathscr{A}$ the unique chamber containing the vertex $v_0$.
Once this chamber has been fixed, we know by \fullref{lem:LinkRootBases}(i) that the sets $\Phi(\mathscr{E})$
and $\Phi(\mathscr{A})$ are in canonical $W$--equivariant bijection. In order to simplify notation, we omit to
write the function $\zeta_C$ which realizes this bijection and identify thereby the sets $\Phi(\mathscr{E})$ and
$\Phi(\mathscr{A})$. In this way, the set $\Pi$ of the root basis $\mathscr{E}$ is identified with $\Psi$, the
set $\Phi(\mathscr{E})^+$ is identified with the set of half-spaces containing $v_0$ and the set
$\Phi(\mathscr{E})^+_{\min}$ with the set of those half-spaces $h$ which contain $v_0$ and such that no wall
separates $v_0$ from $\partial h$.

 We claim that $(\phi_0, \phi_l) = 1$. In order to prove the claim, we will apply
\fullref{prop:FiniteChain}. It follows from the above that $- \phi_l \in \Phi(\mathscr{E})^+_{\min}$.
Write $- \phi_l = \sum_{\psi \in \Psi} \lambda_\psi \psi$ with $\lambda_\psi \geq 0$. Note that by the
definition of $\Psi_0$, we have $\phi_0 \subset \psi$, whence $(\phi_0, \psi) \geq 1$ by
\fullref{lem:LinkRootBases}(iii), for all $\psi \in \Psi \backslash \Psi_0$. Moreover, we have seen that
$\partial \phi_0$ meets $\partial \psi$, whence $|(\phi_0, \psi)| \leq \kappa$, for each $\psi \in \Psi_0$.
Therefore, we have:
$$\begin{array}{rcl}
(\phi_0, \phi_l) &=& (-\phi_0, - \phi_l)\\
 &=& \sum_{\psi \in \Psi} \lambda_\psi (-\phi_0, \psi)\\
 & \leq &  \sum_{\psi \in \Psi_0} \lambda_\psi (-\phi_0, \psi)\\
 & \leq &  \sum_{\psi \in \Psi_0}\kappa \lambda_{\max} \\
 & = & | \Psi_0 | \kappa \lambda_{\max} \\
 & \leq & r \kappa \lambda_{\max}.
 \end{array}$$
Since $l \geq K$, we have $r_l \geq r_K >  r \kappa \lambda_{\max}$ and, therefore, we deduce from
\fullref{thm:chains} that $(\phi_0, \phi_l) = 1$. Moreover, the group $\la r_{\phi_i} \; | \; i=0, \dots,
l\ra$ is infinite dihedral and its parabolic closure $P$ is of irreducible affine type. By \fullref{lem:Pc},
we have $P = \Pc(W(M))$.

It remains to show that $P$ is of rank at least~$3$. By Lemmas~\ref{lem:LinkRootBases}(iv) and
\ref{lem:projection}(i), there exists $\psi \in \Psi$ such that $(\psi, \phi_0) < 0$. We have seen above that if
$\psi \not \in \Psi_0$ then $(\psi, \phi_0) \geq 1$. Thus there exists $\psi \in \Psi_0$ such that $(\psi,
\phi_0) < 0$. We have seen that $\partial \psi$ meets $\partial \phi_i$ for each $i = 0, 1, \dots, l$. Since $l
\geq 8$, we deduce from \fullref{lem:Rk3} that either the reflection $r_\psi$ centralizes $\la r_{\phi_i} \; |
\; i= 0, 1, \dots, l \ra$ or $\la r_\psi, r_{\phi_i} \; | \; i= 0, 1, \dots, l \ra$ is a Euclidean triangle
subgroup. But $r_\psi$ does not commute with $r_{\phi_0}$. Thus by \fullref{lem:Pc} $P$ contains a Euclidean
triangle subgroup and, hence, it is of rank~$\geq 3$.
\end{proof}
\vspace{2pt}

\begin{remark}
Note that a cubical chamber contains finitely many vertices if and only if the $W$--action on $\mathscr{X}$ is
co-compact. Thus \fullref{prop:CubicalChamber} implies that if the $W$--action is not co-compact then $W$
possesses a parabolic subgroup of irreducible affine type and rank~$\geq 3$. Conversely, if $W$ has such a
parabolic subgroup, then it is easily seen that some, whence any, cubical chamber contains infinitely many
vertices. Therefore, we recover the characterization of all those Coxeter groups acting co-compactly on
$\mathscr{X}$; this was first established in \cite{CM05}.
\end{remark}
\vspace{2pt}

Combining the preceding proposition with \fullref{lem:Pc:standard}, one obtains the following useful
precision:
\vspace{2pt}

\begin{cor}\label{cor:Pstandard}
Let $v_0$ be a vertex of $\mathscr{X}_0$ and let $x$ be a vertex of $\mathscr{X}$ belonging to the cubical
chamber containing $v_0$. Let $M \subset \mathscr{M}(v_0, x)$ be a set of pairwise parallel walls of cardinality
greater than $K+8$, where $K$ is the constant of \fullref{prop:CubicalChamber}. Then the parabolic
subgroup $\Pc(W(M))$ is standard with respect to $v_0$.
\end{cor}
\begin{proof}
Up to enlarging $M$ is necessary, we may -- and shall -- assume that $M$ is a maximal subset of
$\mathscr{M}(v_0, x)$ consisting of pairwise parallel walls. By \fullref{prop:CubicalChamber}, the group
$W(M)$ is infinite dihedral and its parabolic closure $P$ is of irreducible affine type (note that enlarging $M$
does not change $P$).

Consider a wall $m$ which separates $v_0$ from some $m' \in M$. Let $M_0$ be the subset of all those elements of
$M_0$ which meet $m$. By \fullref{lem:Rk3}, the set $M_0$ has at most $7$ elements. Therefore, the set $M
\backslash M_0 \cup \{m\}$ is a set of pairwise parallel walls contained in $\mathscr{M}(v_0, x)$, to which we
may apply \fullref{prop:CubicalChamber}. In particular the group $W(M \backslash M_0 \cup \{m\})$ is
infinite dihedral.

Since the set $M \backslash M_0$ contains at least two elements, it follows that $W(M \backslash M_0)$ is
infinite. Therefore, we have $\Pc(W(M \backslash M_0)) = P$ because $W(M \backslash M_0) \subset P$ and any
proper parabolic subgroup of a parabolic subgroup of irreducible affine type is finite. In view of the preceding
paragraph, we deduce from \fullref{lem:Pc}(2) that $r_m \in P$. Since $P$ is of affine type and since $m$ is
parallel to some element of $M$, it follows that $m$ is parallel to all elements of $M$. By the maximality of
$M$, this yields $m \in M$. The desired assertion now follows from \fullref{lem:Pc:standard}.
\end{proof}

Knowing that a parabolic subgroup is standard with respect to $v_0$ will be relevant for the following reason:

\begin{lem}\label{lem:StandardCube}
Let $P$ be a parabolic subgroup of $W$, which is standard with respect to some vertex $v_0$ of the Cayley graph.
Let $M_0$ be a set of walls such that $W(M_0)$ is a finite subgroup of $P$. Then there exists a cube $c_0
\subset \mathscr{X}$ and an element $w \in P$ such that $w.c_0$ contains $v_0$ and that $M(c_0)=M_0$.
\end{lem}
\begin{proof}
Since $W(M_0)$ is finite, so is its parabolic closure $P_0 = \Pc(W(M_0))$. Clearly, it suffices to prove the
lemma for the set $M_0$ consisting of all those walls $m$ such that $r_m \in P_0$.

Since $P_0 \subset P$ and $P$ is standard with respect to $v_0$, there exists $w \in P$ such that $w P_0 w\inv$
is standard with respect to $v_0$. Let $M$ be the set of those walls $m$ such that $r_m \in w P_0w\inv$. Since
$P_0$ is finite, the elements of $M$ meet pairwise by \fullref{lem:PairsHalfspaces}. For each $n \in \N$, let
$M_n$ be the subset of $M$ consisting of the walls at combinatorial distance $n$ from $v_0$ in the Cayley graph
$\mathscr{X}_0$; by convention a vertex is at distance $1$ from a wall if that wall separates the vertex from
one of its neighbors. Note that $W(M_1) = wP_0w\inv$ because $wP_0 w\inv$ is standard with respect to $v_0$. In
particular $M_1$ is nonempty. By \fullref{lem:Pairwise}, there exists a cube $c_1$ containing $v_0$ such that
$M(c_1) = M_1$. Applying \fullref{lem:Pairwise} inductively, one obtains a nested sequence of cubes $c_1
\subset c_2 \subset \dots$ such that $M(c_n) = \bigcup_{i=1}^n M_i$. Since $M$ is finite, the set $M_n$ is empty
for $n$ large enough. Therefore the union $c'_0 = \bigcup_{n} c_n$ is a cube containing $v_0$ and such that
$M(c'_0) = M$. Now the cube $c_0 = w\inv.c'_0$ has the desired property.
\end{proof}

\subsection{Tuples of walls which meet far away from the Cayley graph}

\fullref{thm:FiniteConjug:2sph} of the introduction will be deduced from the following result:

\begin{thm}\label{thm:cubes}
There exists a constant $A=A(W,S)$ such that the following property holds. Let $M$ be a set of walls such that
the intersection $\bigcap_{m \in M} m$ is nonempty in $\mathscr{X}$.
If the distance from $\mathscr{X}_0$ to $\bigcap_{m \in M} m$ is at least  $A$,
then $W(M)$
has a direct component of affine type and rank~$\geq 3$; in particular, it has a free abelian normal subgroup of
rank $\geq 2$.
\end{thm}
\begin{proof}
By \fullref{lem:NR} combined with Ramsey's theorem, there exists a constant $K'$ such that any set of at least
$K'$ walls contains a subset of more than $K + 8$ pairwise parallel walls, where $K$ is the constant of
\fullref{prop:CubicalChamber}. We choose $A \in \R_+$ large enough so that the ball of combinatorial
radius $K'$ centered at some vertex of $\mathscr{X}_0$ is properly contained in the ball of radius $A$ centered
at that same vertex. Since $W$ acts transitively on the vertices of $\mathscr{X}_0$, the so-defined constant $A$
does not depend on the chosen vertex.

By \fullref{lem:Pairwise} there exists a cube $c \subset \mathscr{X}$ such that $M(c) = M$. In order to prove
the theorem, it suffices to show that if every such cube is at combinatorial distance at least $K'$ from
$\mathscr{X}_0$, then $W(M)$ has a direct component of affine type and rank~$\geq 3$.

Let thus $c$ be a cube such that $M(c) = M$ and that $c$ is at minimal combinatorial distance from
$\mathscr{X}_0$. Let $d$ be this distance, let $x$ be a vertex of $c$ and $v_0$ be a vertex of $\mathscr{X}_0$
at combinatorial distance $d$ from $x$ and assume that $d> K'$.

Let $\Psi$ be the set of those half-spaces which contain $v_0$ but not its neighbors in $\mathscr{X}_0$. Assume
that that $x \not \in \psi$ for some  $\psi \in \Psi$. Given any minimal path from $x$ to $v_0$, this path
crosses the wall $\partial \psi$. Since $\partial \psi$ separates $v_0$ from one of its neighbors, say $v$, it
follows that there exists a minimal path from $x$ to $v_0$ whose last edge crosses $\partial \psi$. Thus the
distance from $x$ to $v$ is one less than the distance from $x$ to $v_0$. This contradicts the definition of
$v_0$ since, by the definition of the elements of $\Psi$, the vertex $v$ belongs to $\mathscr{X}_0$. This shows
that $x$ belongs to the cubical chamber containing $v_0$.

The cardinality of $\mathscr{M}(x, v_0)$ coincides with the combinatorial distance from $x$ to $v_0$. By the
definition of $K'$, it follows that $\mathscr{M}(x, v_0)$ contains a subset $M'$ of pairwise parallel walls and
of cardinality greater than $K+8$. Up to enlarging $M'$ is necessary, we may -- and shall -- assume that $M'$ is
a maximal subset of pairwise parallel walls contained in $\mathscr{M}(v_0, x)$. By
\fullref{prop:CubicalChamber}, the group $W(M')$ is infinite dihedral and its parabolic closure, which
we denote by $P$, is of irreducible affine type.

Let $x'$ be a vertex of the cube  $c$ neighboring $x$ and let $m \in M$ be the wall which separates $x$ from
$x'$. By definition, the combinatorial distance from $x'$ to $\mathscr{X}_0$ is at least $d$. This implies that
the combinatorial distance from $v_0$ to $x'$ is $d+1$ and, hence, that $\mathscr{M}(v_0,x') = \mathscr{M}(v_0,
x) \cup \{m\}$. We now show that the reflection $r_m$ normalizes the parabolic subgroup $P$.

Assume first that $m = \partial \psi$ for some $\psi \in \Psi$. In that case, the wall $m$ must meet every
element of $M'$ otherwise $m$ would be separated from $\partial \psi$ by some element of $M'$, which is absurd
since $m = \partial \psi$. By \fullref{lem:Rk3}, this implies that either $r_m$ centralizes $W(M')$ or that
$W(M' \cup \{m\})$ is a Euclidean triangle subgroup. In the first case, the reflection $r_m$ normalizes $P$; in
the second one, we have $r_m \in P$ by \fullref{lem:Pc}.

Assume now that $m \not \in \{\partial \psi \; | \; \psi \in \Psi\}$. Equivalently, this means that $x'$ belongs
to the cubical chamber containing $v_0$.

Suppose that $m$ does not meet any element of $M'$. Then $M' \cup \{m\}$ is a set of pairwise parallel walls to
which we may apply \fullref{prop:CubicalChamber}. This proves that $W(M' \cup \{m\})$ is infinite
dihedral and, hence, the reflection $r_{m}$ through $m$ belongs to $P$ by \fullref{lem:Pc}.

Suppose now that the subset $M'(m)$ of those elements of $M'$ which meet $m$ is nonempty. Suppose first that
$M'(m)$ contains less than $8$ elements. Then  
$$M' \backslash \big(M'(m)\big)\cup \{m\} \subset
\mathscr{M}(v_0, x')$$ 
is a set of pairwise parallel walls to which we may apply
\fullref{prop:CubicalChamber}. As in the preceding paragraph, we obtain that $r_m$ belongs to $P$.
Suppose now that $M'(m)$ contains at least $8$ elements. By \fullref{lem:Rk3}, this implies that either $r_m$
centralizes $W(M')$ or that $W(M' \cup \{m\})$ is a Euclidean triangle subgroup. In the first case, the
reflection $r_m$ normalizes $P$; in the second one, we have $r_m \in P$ by \fullref{lem:Pc}.

In all cases, we have seen that $r_m$ normalizes $P$. Since every element of $M$ separates $x$ from one of its
neighboring vertices in $c$, it follows that the group $W(M)$ is contained in the normalizer $N_W(P)$ of $P$ in
$W$.
Therefore, by \fullref{lem:Normalizer}, each direct component of $W(M)$ is either contained in $P$ or
centralizes $P$.

Assume that every wall of $M$ meets every wall of $M'$. By the maximality of $M'$, the set $M'$ possesses an
element $m'$ such that $m'$ is not separated from $x$ by any wall. By the definition of $\mathscr{X}$, it
follows that $x$ belongs to an edge which is cut by the wall $m'$. Since every wall in $M$ meets $m'$, we deduce
from \fullref{lem:Pairwise} that $c$ is a face of a $(n+1)$--cube $c''$ whose center is contained in $m'$. Let
$c'$ be the $n$--cube which is opposite $c$ in $c''$. Thus $c$ and $c'$ are separated by $m'$ and $M(c) = M(c') =
M$. By construction, the combinatorial distance from $v_0$ to $c'$ is strictly smaller than the combinatorial
distance from $v_0$ to $c$, which contradicts the definition of $c$.

Therefore $M$ contains an element $m_P$ which does not meet all elements of $M'$. By the same arguments as
above, we see that $r_{m_P}$ belongs to $P$. Since $P$ is a Coxeter group of affine type in which the set $M'$
corresponds to one direction of hyperplanes, it follows that $m_P$ does not meet any element of $M'$. If $m \in
M$ is another element of $M$ which does not meet all elements of $M'$, then we obtain similarly that $r_m \in P$
and that $m$ is parallel to all elements of $M' \cup \{m_P\}$. Since $m$ meets $m_P$ because the elements of $M$
meet pairwise, it follows that $m= m_P$. This proves that, with the exception of $m_P$ which meets no element of
$M'$, all other elements of $M$ meet all elements of $M'$.

Let $P_0$ be the subgroup generated by the direct components of $W(M)$ which are contained in $P$. Thus $r_{m_P}
\in P_0 \subset P$. Note that any reflection subgroup of an affine Coxeter group is either finite or of affine
type. Therefore, in order to finish the proof, it suffices to prove that $P_0$ is infinite because $P_0$ is of
$2$--spherical type by \fullref{lem:W(M):2sph} and, hence, if it is infinite, then its rank is at least~$3$.

Let $M_0 = \{ m \in M \; | \; r_m \in P_0\}$ and let $M_1 = M \backslash M_0$. Thus $W(M_1)$ centralizes $P$ and
every element of $M_1$ meets every element of $M'$. Recall that $d$ denotes the combinatorial distance from
$v_0$ to $x$. By \fullref{lem:Pairwise} there is a cube $c'_1$, containing $x$, such that $M(c'_1) = M_1 \cup
\{m'\}$. Let $c_1$ be the face of that cube such that $x \not \in c_1$ and $M(c_1) = M_1$. In particular, the
combinatorial distance from $c_1$ to $v_0$ equals $d-1$.

Assume now that $P_0 = W(M_0)$ is finite in order to obtain a contradiction. By \fullref{lem:StandardCube},
there exists an element $w \in P$ and a cube $c_0 \subset \mathscr{X}$ such that $M(c_0) = M_0$ and that $w
.c_0$ contains $v_0$. Note that $M(w .c_0 ) = w(M(c_0))$ and hence $W(M(w.c_0))=w W(M_0) w \inv \subset P$.
Since $P$, and hence $w$, centralizes $W(M_1)$ it follows that every wall in $M_1$ meets every wall in
$M(w.c_0)$. Note that if $c, c'$ are cubes of $\mathscr{X}$ at combinatorial distance $r$ from one another and
such that $M(c) \subset M(c')$, then every vertex of $c$ is at combinatorial distance $r$ from $c'$. Therefore,
by \fullref{lem:Pairwise}, there is a cube $c_2$, at combinatorial distance at most $d-1$ from $v_0$, such
that $M(c_2) = M(w.c_0) \cup M(c_1)$. Since $w$ centralizes $W(M_1) = W(M(c_1))$, if follows that $w(M(c_1)) =
M(c_1)$.  Therefore, the cube $w\inv .c_2$ is at combinatorial at most $d-1$ from the Cayley graph
$\mathscr{X}_0$ and we have $M(w\inv .c_2) = M = M(c)$. This contradicts the definition of $c$, which finishes
the proof.
\end{proof}

\subsection[Proof of \ref{thm:FiniteConjug:2sph}]{Proof of \fullref{thm:FiniteConjug:2sph}}

Let $\mathscr{T}(W)$ be the set of all subsets $T \subset W \backslash \{1\}$ such that each pair of elements of
$T$ generates a finite group. Let also $\mathscr{C(X)}$ denote the union of all cubes of $\mathscr{X}$. We
define a map $\sigma \co  \mathscr{T}(W) \to \mathscr{C(X)}$ as follows.

Let $t \in T$. Thus $t$ is of finite order and, hence, its parabolic closure $\Pc(t)$ is finite. Define $M(t) :=
\{ m \in \mathscr{M(A)} \; | \; r_m \in \Pc(t) \}$. Thus for all $m, m' \in M(t)$, the group $\la r_m, r_m' \ra$
is finite and, hence, the wall $m$ meets $m'$. Let now $t, t' \in T$. Since $\la t,t' \ra$ is finite, so is its
parabolic closure $\Pc(\la t, t' \ra)$, which contains $\Pc(t)$ and $\Pc(t')$ by definition. Therefore, for any
$m \in M(t)$ and $m' \in M(t')$, the group $\la r_m, r_{m'} \ra$ is finite and, hence, the wall $m$ meets $m'$.
This shows that the elements of $M(T) = \bigcup_{t \in T} M(t)$ meet pairwise.

By \fullref{lem:Pairwise}, there is a cube $c$ in $\mathscr{X}$ such that $M(c) = M(T)$. Among all such cubes,
choose one which is at minimal combinatorial distance from the Cayley graph $\mathscr{X}_0$; we define
$\sigma(T)$ to be that cube. Note that the group $W$ acts on $\mathscr{T}(W)$ by conjugation and on
$\mathscr{C(X)}$ via its action on $\mathscr{X}$, but the map $\sigma$ is \emph{not} $W$--equivariant because the
definition of $\sigma$ depends on an arbitrary choice.

Now we associate to each cube $c \in \mathscr{C(X)}$ a finite subset $\mathscr{T}(c) \subset \mathscr{T}(W)$ as
follows. Consider a $k$--tuple of subsets $M_1, M_2, \dots M_k \subset M(c)$ which satisfy the following
conditions:
\begin{itemize}
\item $M(c) = \bigcup_{i=1}^k M_i$;

\item For all $i, j \in \{1, \dots, k\}$, the group $W(M_i \cup M_j)$ is finite.
\end{itemize}

Given such a $k$--tuple, choose a nontrivial element $t_i \in W(M_i)$ for each $i= 1, \dots, k$. Clearly we have
$T = \{t_i \; | \; i=1, \dots, k\} \in \mathscr{T}(W)$. We denote by $\mathscr{T}(c)$ the set consisting of all
those elements of $\mathscr{T}(W)$ which are obtained from $c$ in this manner. Note that the construction of $T$
depends on some choices, but each choice has to be made between a finite number of possibilities. Therefore
$\mathscr{T}(c)$ is a finite subset of $\mathscr{T}(W)$. Note also that the map $\mathscr{T} \co \mathscr{C(X)}
\to 2^{\mathscr{T}(W)}$ is $W$--equivariant: for all $w \in W$ we have $\mathscr{T}(w.c) = \{ w T w\inv \; | \; T
\in \mathscr{T}(c)\}$.

Let $\mathscr{C}_0(\mathscr{X})$ be the set of all those cubes $c$ such that $W(M(c))$ has no direct component
of affine type and that $c$ is at minimal combinatorial distance from $\mathscr{X}_0$ among all cubes $c'$ such
that $M(c') = M(c)$. Clearly the $W$--action on $\mathscr{X}$ preserves $\mathscr{C}_0(\mathscr{X})$. By
\fullref{thm:cubes}, the distance from $\mathscr{X}_0$ to any element of $\mathscr{C}_0(\mathscr{X})$ is
uniformly bounded. Since $\mathscr{X}$ is locally finite and $W$ is transitive on the vertices of
$\mathscr{X}_0$, it follows that $W$ has finitely many orbits in $\mathscr{C}_0(\mathscr{X})$.

Let now $\mathscr{T}_0(W) \subset \mathscr{T}(W)$ be the subset consisting of all those $T$ such that $W(M(T))$
has no direct component of affine type. Note that $W$ acts on $\mathscr{T}_0(W)$ by conjugation. We claim that
$W$ has finitely many orbits in $\mathscr{T}_0(W)$. Let $\{c_1, \dots, c_k \} \subset
\mathscr{C}_0(\mathscr{X})$ be a set of representatives of the $W$--orbits and let $T \in \mathscr{T}(W)$. By
definition, we have $\sigma(T) \in \mathscr{C}_0(\mathscr{X})$ and $T \in \mathscr{T}(\sigma(T))$. Let $w \in W$
such that $w.\sigma(T) = c_i$ for some $i$. Thus $wTw\inv \in \mathscr{T}(c_i)$. This shows that the finite
subset $\bigcup_{i =1} ^k \mathscr{T}(c_i) \subset \mathscr{T}(W)$ contains a representative of each $W$--orbit
in $\mathscr{T}_0(W)$. This proves the claim.

 Let $\mathscr{G}$ be any set of subgroups of $W$ invariant under conjugation. Assume that each
$\Gamma \in \mathscr{G}$ possesses a generating set which belongs to $\mathscr{T}_0(W)$. Then, by the claim, the
set $\mathscr{G}$ is a \emph{finite} union of conjugacy classes. Using this observation, we now prove the
desired assertions successively.

(i)\qua The `only if' part is clear. Suppose that $W$ has a parabolic subgroup of irreducible affine type
and rank~$\geq 3$. Then $W$ has no reflection subgroup of irreducible affine type and rank~$\geq 3$ by
\fullref{prop:PcAffine} and it follows that $\mathscr{T}_0(W) = \mathscr{T}(W)$. By the above, it
follows that $W$ has finitely many conjugacy classes of  $2$--spherical subgroups.

(ii)\qua Let $\mathscr{G}_1$ be the set of all $2$--spherical reflection subgroups with no direct
component of irreducible affine type. By definition, every $\Gamma \in \mathscr{G}_1$ has a generating set $T
\in \mathscr{T}(W)$ consisting of reflections, and such that $\Gamma = W(M(T))$. Therefore $T \in
\mathscr{T}_0(W)$ and the desired finiteness property follows.

(iii)\qua Let $\mathscr{G}_2$ be the set of all $2$--spherical subgroups $\Gamma$ such that $\Gamma$ has no
nontrivial free abelian normal subgroup and $\Zrk(\Gamma) \geq \Zrk(W)-1$. Let $\Gamma \in \mathscr{G}_2$ and $T
\in \mathscr{T}(W)$ be a generating set of $\Gamma$. Clearly $\Gamma \subset W(M(T))$. Let $W(M(T))=W_1 \times
\dots \times W_l$ be the decomposition of $W(M(T))$ into its direct components, and assume that $W_1$ is of
affine type. Let $T_1$ be the translation subgroup $W_1$. We have $T_1 \cap \Gamma = \{1\}$ because $T_1 \cap
\Gamma$ is a free abelian normal subgroup of $\Gamma$. On the other hand, it follows from the direct
decomposition above that any free abelian subgroup of $\Gamma$ has a finite index subgroup which centralizes
$T_1$. Moreover $W(M(T))$ is $2$--spherical by \fullref{lem:W(M):2sph}; therefore, as an infinite $2$--spherical
Coxeter group, $W_1$ is of rank~$\geq 3$ and hence $\Zrk(T_1)\geq 2$. This shows that $\Zrk(\la T_1 \cup \Gamma
\ra) \geq \Zrk(\Gamma) + 2
> \Zrk(W)$, a contradiction. Therefore $W(M(T))$ has no direct component of affine type. In other words,
we have $T \in \mathscr{T}_0(W)$ and the desired finiteness property
follows.

(iv)\qua Let $\mathscr{G}_3$ be the set of all $2$--spherical subgroups $\Gamma$ such that $\Gamma$ is not infinite
virtually abelian. Let $\Gamma \in \mathscr{G}_3$ and $T \in \mathscr{T}(W)$ be a generating set of $\Gamma$.
Let $W(M(T))=W_1 \times \dots \times W_l$ be the decomposition of $W(M(T))$ into its direct components, and
assume that $W_i$ is of affine type for each $i \leq j$ but $W_i$ is not for $i > j$. Thus for all $i \leq j <
i'$, the group $W_{i'}$ centralizes $W_i$ and, hence, normalizes its parabolic closure $\Pc(W_i)$, which is of
affine type by \fullref{prop:PcAffine}. Therefore, by \fullref{lem:Normalizer}, either $W_{i'}$
centralizes $\Pc(W_i)$ or $W_{i'} \subset \Pc(W_i)$. Since the centralizer of $\Pc(W_i)$ is finite by hypothesis
and since $W_{i'}$ is not affine, it follows in both cases that $W_{i'}$ is finite. This shows that if $W(M(T))$
has a component of affine type, then each direct component of $W(M(T))$ is either affine or finite. In that case
$W(M(T))$ would be virtually abelian, which is impossible since $\Gamma \subset W(M(T))$. Therefore $W(M(T))$
has no component of affine type and $T \in \mathscr{T}_0(W)$. \qed

\subsection[Co-Hopfian Coxeter groups : proof of \ref{cor:Hopfian}]{Co-Hopfian Coxeter groups : proof of \fullref{cor:Hopfian}}

Let $W = W_1 \times \dots \times W_l$ be the decomposition of $W$ into its direct components.

If $W_1$ is of affine type, then there exists a monomorphism $\phi_1 \co  W_1 \to W_1$ which is not surjective.
Then the unique homomorphism $\phi$ whose restriction on $W_1$ (resp.\ $W_i$) is $\phi_1$ (resp.\ the identity for
$i > 1$) is a monomorphism which is not surjective. Thus $W$ is not co-Hopfian.

Assume now that for each $i$ the group $W_i$ is not affine. By \fullref{lem:FreeAbNormal}, this implies that
$W$ has no nontrivial free abelian normal subgroup. Therefore, the group $W$ contains finitely many conjugacy
classes of subgroups isomorphic to $W$ by \fullref{thm:FiniteConjug:2sph}(iii). By the main result of
\cite{HRT97}, the outer automorphism group of $W$ is finite. We deduce that $W$ admits only finitely many
monomorphisms into itself up to conjugation. The rest of the proof is similar to \cite[Proof of
Theorem~3.1]{RS94}; for convenience, we reproduce the details.

Let $\varphi \co  W \to W$ be a monomorphism and assume that $\varphi$ is not surjective in order to obtain a
contradiction. For each $n$, the centralizer of $\varphi^n(W)$ in $W$ is finite, otherwise it would contain an
element of infinite order $\gamma$ (by Selberg's lemma, $W$ has a torsion free subgroup of finite index) and the
group $\la \{\gamma\} \cup \varphi^n(W) \ra$ would be of $\Z$--rank strictly greater than $\Zrk(W)$. Since $W$
has finitely many conjugacy classes of finite subgroups, there is a finite subgroup $A$ which is the centralizer
of $\varphi^n(W)$ for all sufficiently large $n$. In particular $\varphi(A) = A$.

Since $W$ admits only finitely many monomorphisms into itself up to conjugation, there exist arbitrarily large
integers $k, l$ and an element $g \in W$ such that for all $w\in W$ we have
$$g \varphi^k(w) g\inv = \varphi^{k+l}(w).$$
Therefore, conjugating $\varphi^k(w)$ by $g$ is equivalent to transform it by $\varphi^l$, so we have:
$$g \varphi^{k+l}(w) g\inv = \varphi^l(\varphi^{k+l}(w)) = \varphi^l(g \varphi^k(w) g\inv ) = \varphi^l(g)
\varphi^{k+l}(w) \varphi^l(g\inv).$$ In particular $g\inv \varphi^l(g) \in A$, so $\varphi^l(g) = ga$ for some
$a \in A$. Since $\varphi(A) = A$, for each $m \geq 1$ we have $\varphi^{ml}(g) = g a_m$ for some $a_m \in A$.
By the pigeonhole principle, we obtain $\varphi^{m_1 l}(g) = \varphi^{m_2 l}(g)$ for some $m_1 < m_2$. But this
shows that $\varphi^{m_1 l}(g) \in \varphi^n(W)$ for all integers $n$. As for all $w \in W$ we have:
$$\varphi^{k+l+m_1 l}(w) =\varphi^{m_1 l}(g) \varphi^{k + m_1 l}(w) \varphi^{m_1 l}(g\inv)$$
and since moreover $\varphi^{m_1 l}(g) \in \varphi^{k + m_1 l}(W)$, it follows that $$\varphi^{k+l+m_1 l}(W)
=\varphi^{k + m_1 l}(W),$$ a contradiction. \qed

\medskip 
\textbf{Acknowledgement}\qua The author is supported by a fellowship
of `Aspirant' from the F.N.R.S.  (Fonds National de la Recherche
Scientifique -- Belgium).

\bibliographystyle{gtart}
\bibliography{link}

\end{document}